
\input amstex
\baselineskip=18pt
\magnification\magstephalf
\parskip=8pt plus 1pt minus 1pt
\TagsOnRight

\def \Ai{\text{\rm Ai}\,}
\def \Vol{\,\text{\rm Vol}\,}
\def \supp{\,\text{\rm supp}\,}
\def \grad{\,\text{\rm grad}\,}
\def \sgn{\,\text{\rm sgn}\,}
\def \rat{\text{\rm rat}}

\def \const{\,\text{\rm const}\,}
        
        \def \dist{\,\text{\rm dist}\,}
        \def \Area{\,\text{\rm Area}\,}

        \def \CS{\cos(2\pi RY_0(n)-\phi)}
        \def \CSpm{\cos(2\pi RY_\pm(n)-\phi_\pm)}

        \def \PI{\pi^{-1}}

        \def \II{\int_{-\infty}^\infty}

        \def \sz{\sum_{n\in\Z^2\setminus\{0\}}}

        \def \G{\Gamma}
        \def \g{\gamma}
\def \a{\alpha}
        \def \b{\beta}
        \def \de{\delta}
        \def \De{\Delta}
        \def \ep{\varepsilon}
        \def \la{\lambda}
        \def \La{\Lambda}
        
        \def \t{\theta}
        
        \def \sg{\sigma}
        \def \Om{\Omega}
        \def \om{\omega}

        \def \f{\varphi}
        
        \def \Q{\bold Q}
        \def \R{\bold R}

        \def \Z{\bold Z}

        \def \lT{\limsup_{T\to\infty}}
        \def \liT{\lim_{T\to\infty}}

        \def \lL0{\lim_{S\to\infty,\,S/T\to 0}}
        \def \lt0{\lim_{T\to\infty,\,S/T^2\to 0}}

        \def \di{\displaystyle}
        \def \IT0{{1\over T}\int_{T_0}^T}
        \def \IT{{1\over T}\int_1^T}
        \def \IO{{1\over T}\int_0^T}

        \def \tc{\tilde\chi}
        \def \tQ{\tilde Q}
        \def \tp{\tilde\varphi}
        
        \def \tpd{\tilde\varphi(2\pi n\de)}
        \def \pn{\psi_l(n)}
        \def \sZ{\sum_{n\in\bold Z^2}}

\def \ex{\exists\,}
\def \SL{\sum_{n\in\Z^2\setminus L_\g}}

\def \ar{\overline}
\def \RHS{\,\text{\rm RHS}\,}
\def \Gjl{\G_j^{(l)}}
\def \Elm{E_{lm}(R)}
\def \Flm{F_{lm}(R)}
\def \Plm{\Phi_{lm}(R)}
\def \np{n_2^\perp}
\def \xp{\xi_2^\perp}
\def \BS{\text{\rm BS}}

\null
 \rightline {IASSNS-HEP-93/35}

\rightline {May 1993}

\bf
\centerline {Distribution of Energy Levels}

\vskip 3mm

\centerline{of a Quantum Free Particle on a}

\vskip 3mm

\centerline{Surface of Revolution}

\vskip 5mm

\centerline {Pavel M. Bleher}

\vskip 3mm

\centerline {School of Natural Science}

\centerline {Institute for Advanced Study}

\centerline {Princeton, NJ 08540}

\rm

\vskip 1cm

{\bf Abstract}. We prove that the error term $\La(R)$ in the Weyl
asymptotic formula
$$
\#\{ E_n\le R^2\}={\Vol M\over 4\pi} R^2+\La(R),
$$
for the Laplace operator on a surface of revolution $M$
satisfying a twist hypothesis, has the form $\La(R)
=R^{1/2}F(R)$ where $F(R)$ is an almost periodic
function of the Besicovitch class $B^2$, and the Fourier
series of $F(R)$ in $B^2$ is $\sum_\g A(\g)\cos(|\g|R-\phi)$
where the sum goes over all closed geodesics on $M$, and
$A(\g)$ is computed through simple geometric characteristics
of $\g$. We extend this result to surfaces of revolution,
which violate the twist hypothesis and satisfy a more general
Diophantine hypothesis. In this case we prove that
$\La(R)=R^{2/3}\Phi(R)+R^{1/2}F(R)$, where $\Phi(R)$ is a
finite sum of periodic functions and $F(R)$ is an almost periodic
function of the Besicovitch class $B^2$. The Fourier
series of $\Phi(R)$ and $F(R)$ are computed.

\vfill\eject

\beginsection 1. Introduction \par

Let $M$ be a two--dimensional smooth compact
manifold which is homeomorphic to a sphere, and which
is a surface of revolution in $\R^3$, with an axis $A$
and poles $N$ and $S$ (see Fig.1). The geodesic flow on $M$
is a classical integrable system due to the Clairaut
integral,
$$
r\sin\a=\const.\eqno (1.1)
$$
In the present work we are interested in high energy levels
of the corresponding quantum system,
$$
-\De u_n= E_nu_n.\eqno (1.2)
$$
Let $s$ be the normal coordinate (the length of geodesic)
along meridian and
$$
r=f(s),\quad 0\le s\le L,\eqno(1.3)
$$
be the equation of $M$, where $r$ is the radial coordinate.
Then
$$
\De=f(s)^{-1}{\partial\over\partial s}
\left( f(s){\partial\over\partial s}\right)
+f(s)^{-2}{\partial^2\over\partial\f^2},
\eqno (1.4)
$$
where $\f$ is the angular coordinate.

We will assume that $f(s)$ has a simple structure, so that
$$
f'(s)\not=0,\quad s\not= s_{\max};\quad\quad
f''(s_{\max})\not=0,\eqno (1.5)
$$
where
$$
f(s_{\max})=\max_{0\le s\le L}f(s)\equiv f_{\max}.
$$
For normalization we put $f_{\max}=1.$
Another assumption on $M$ is the following twist hypothesis.

Consider the equator on $M$,
$$
\g_E=\{ s=s_{\max},\;\; 0\le\f\le 2\pi\},
$$
(we do not assume that $M$ is symmetric with respect to
$\g_E$, but we still call $\g_E$ the equator
keeping a visual interpretation of objects on $M$),
and a geodesic $\g$ which starts at
$x_0=(s=s_{\max},\f=0)\in\g_E$ at some angle
$-\pi/2<\a_0<\pi/2$ to the direction to the north.
The Clairaut integral on $\g$ is $I=\sin\a_0$,
and we can parametrize $\g$ by $-1< I<1\:\;\;
\g=\g(I)$. It follows from the Clairaut integral
that $\g(I)$ oscillates between two parallels,
$s=s_+$ and $s=s_-$, where $f(s_-)=f(s_+)=I$, so
$\g(I)$ intersects $\g_E$ infinitely many times.
Let $x_n$ be the $n$--th intersection of $\g$
with $\g_E$, $n\in\Z$. Define
$$
\tau(I)=|\g[x_0,x_2]|,
$$
the length of $\g$ between $x_0$ and $x_2$,
and
$$
\om(I)=(2\pi)^{-1}(\f(x_2)-\f(x_0)),
$$
the phase of $\g$ between $x_0$ and $x_2$ (see Fig. 1).
Observe that $\om(I)$ is defined $\mod 1$. To define $\om(I)$
uniquely, we choose a continuous branch of $\om(I)$
starting at $\om(0)=0.$ Then $\om(-I)=-\om(I)$ and for $I\ge 0$,
$$
\om(I)=\pi^{-1}\int_{s_-}^{s_+}{d\f\over ds}\,ds-1.
\eqno (1.6)
$$
Define $\tau(1)=\lim_{I\to 1-0}\tau(1)$ and
$\om(1)=\lim_{I\to 1-0}\om(I).$

It is easy to see that a finite geodesic $\g$
with the Clairaut integral
$0<I<1$ is closed iff $\om(I)$ is rational.
More precisely, let $n(\g)$ denote the number
of revolutions of a closed geodesic
$\g$ around the axis $A$
and $m(\g)$ denote the number of oscillations
of $\g$ along meridian. Then
$$
\om(I)=(n(\g)/m(\g))-1.
$$
To facilitate formulation of subsequent results
we take the convention that a finite geodesic
$\g$ with $I=1$, which goes along the equator, is closed
iff both $n(\g)$ and
$$
m(\g)\equiv {n(\g)\over \om(1)+1}
$$
are integers.

{\bf Twist Hypothesis} (TH). $\om'(I)\not=0$,
$\forall\, I\in[0,1]$.

To illustrate TH consider an ellipsoid of revolution,
$$
{x^2\over a^2}+{y^2\over a^2}+{z^2\over b^2}=1.
$$
If $a<b$ (oblong ellipsoid), then $\om'(I)>0$,
if $a>b$ (oblate ellipsoid), then $\om'(I)<0$;
so TH holds in both cases.
Curves (a) and (b) on Fig.2 are the graphs of $\om(I)$ found with the help of
computer for the ellipsoids of revolution with $a=1$, $b=2$
and $a=2$, $b=1$, respectively.
The cross--sections of the ellipsoids are shown in the lower part
of Fig.2.
TH is violated
for a sphere ($a=b$), when $\om'(I)\equiv 0$.
TH can be also violated
for a bell--like shape of $M$ shown on Fig.2 (see the cross--section
(c) and the graph (c) of $\om(I)$ on this figure) and in some other cases.

Let
$$
N(R)=\#\{ E_n\le R^2\}
\eqno (1.7)
$$
be the counting function of $E_n$. Then the Weyl law
says that
$$
N(R)={\Vol M\over 4\pi}\, R^2+\La(R)
$$
where
$\La(R)=o(R^2),\; R\to\infty$.
A general estimate of H\"ormander [H\"or1]
gives $$
\La(R)=O(R).
$$
This estimate is sharp for $S^2$ and some other degenerate
surfaces for which closed geodesics cover a set of positive
Liouville measure in the phase space. If the Liouville
measure of the union of all closed geodesics in the
phase space is 0, then as was shown by
Duistermaat and Guillemin [DG],
$$
\La(R)=o(R).
$$
For surfaces of negative curvature Selberg and B\'erard [B\'er]
proved a better estimate:
$$
\La(R)=O(R/\log R),
$$
and it is a very difficult open problem to show that
$\La(R)=O(R^{1-\ep})$ for some $\ep>0$, even in the case
of constant negative curvature (see recent works [Sar],
[LS] and [HR] where statistics of eigenvalues
and eigenfunctions of the Laplace operator on surfaces
of constant negative curvature is discussed).

For a flat torus $\La(R)$ reduces to the error term of the classical
circle problem, and the best estimate here is due to Huxley [Hux]:
$$
\La(R)=O(R^{46/73}(\log R)^{315/146}).
$$
A well--known conjecture of Hardy [Har1]
$$
\La(R)=O(R^{(1/2)+\ep}),\quad\forall\ep>0,
$$
is probably also a very difficult open problem. On the
other hand Hardy proved [Har2] that
$$
\limsup_{R\to\infty}R^{-1/2}|\La(R)|=\infty,
$$
so $(1/2)+\ep$ is the best possible exponent.

Colin de Verdi\`ere
[CdV1,CdV2] proved that for a generic surface
of revolution of simple structure
$$
\La(R)=O(R^{2/3}).\eqno (1.8)
$$
We prove in the present paper the following result:

{\bf Theorem 1.1.} {\it Assume that $M$ is a surface of
revolution of simple structure, and $M$ satisfies TH.
Then
$$
N(R)={\Vol M\over 4\pi}\, R^2+
R^{1/2}F(R),\eqno (1.9)
$$
where $F(R)$ is an almost periodic function of
the Besicovitch class $B^2$, and the Fourier
series of $F(R)$ in $B^2$ is
$$
F(R)=\sum_{\text{\rm closed geodesics $\g$}}
A(\g)\cos(|\g|R-\phi),\eqno (1.10)
$$
where summation goes over all closed (in general, multiple)
oriented geodesics $\g\not=0$
on $M$, $\phi=(\pi/2)+(\pi/4)\sgn \om'(I)$,
and
$$\eqalign{
A(\g)
&=\PI(-1)^{m(\g)}|\om'(I)|^{-1/2}m(\g)^{-3/2}\cr
&=\PI(-1)^{m(\g)}|\om'(I)|^{-1/2}\tau(I)^{3/2}
|\g|^{-3/2},\quad I=I(\g).\cr}\eqno (1.11)
$$}

In Theorem 1.2 we extend Theorem 1.1 to the case when TH
is violated. In this case we introduce

{\bf Diophantine Hypothesis} (DH). {\it Assume that $\om(I)$ has
at most finitely many critical points
$0<I_1<\dots<I_K<1$ (so that $I=0,1$ are not critical)
with $\om'(I_k)=0$,
and $\om''(I_k)\not=0$, $k=1,\dots,K$. Assume, in addition, that
for every $k=1,\dots,K$, $\om(I_k)$ is either rational or
Diophantine in the sense that $\ex 1>\zeta>0$ and $C>0$ such that
$$
\left|\om(I_k)-{p\over q}\right|\ge {C\over q^{2+\zeta}},
\quad\forall\;{p\over q}\in\Q.
\eqno (1.12)
$$}

{\bf Theorem 1.2.} {\it Assume that $M$ is a surface of revolution of
simple structure and DH holds.  Then
$$
N(R)={\Vol M\over 4\pi}R^2+
R^{2/3}\sum_{k\: \om(I_k)\in\Q}\Phi_k(R)
+R^{1/2}F(R),\eqno (1.13)
$$
where $\Phi_k(R)$ are bounded periodic functions,
and $F(R)$ is an almost periodic function of
the Besicovitch class $B^2$.
The Fourier series of $\Phi_k(R)$ is
$$
\Phi_k(R)=(1/2)3^{-2/3}\G(2/3)\pi^{-4/3}\tau(I_k)^{4/3}
\sum_{\g\: I(\g)=I_k}
(-1)^{m(\g)}|\g|^{-4/3}\sin(|\g|R),
\eqno (1.14)
$$
and the Fourier series of $F(R)$ is
$$
F(R)=\sum_{\g\: I(\g)\not= I_1,\dots,I_K}
A(\g)\cos(|\g|R-\phi(\g)),\eqno (1.15)
$$
where $\phi(\g)=(\pi/2)+(\pi/4)\sgn \om'(I)$, $I=I(\g)$, and
$A(\g)$ is given in (1.11).}

In the works [H-B], [BCDL], [Ble1], [Ble2] and [BL] some general
results were proved on the existence and properties
of a limit distribution of any almost periodic function
of the Besicovitch class $B^2$ (see especially
Theorems 4.1--4.3 in [Ble1] and Theorems 3.1, 3.3
in [Ble2]). Theorem 1.2 combined with these results
lead us to the following

{\bf Corollary.} {\it Assume that $M$ is a surface of revolution
of simple structure and $M$ satisfies TH or, more
generally, DH.
Then the normalized error function $F(R)$ in the formulas (1.9) and (1.13)
has a limit distribution
$\nu(dt)$, i.e., for every bounded continuous function $g(t)$ on the line,
$$
\liT\IO g(F(R))\,dR=\int_{-\infty}^\infty g(t)\,\nu(dt).
$$
If, in addition,
the lengths of all primitive closed geodesics
on $M$ with $I\ge 0$ are linearly independent over $\Z$,
then $\nu(dt)$ is absolutely continuous and the
density function $p(t)=\nu(dt)/dt$ is an entire function of $t$
which satisfies on the real axis the estimates
$$
\eqalign{
&0\le p(t)\le C\exp(-\la t^4),\quad C,\la>0,\cr
&P(-t),\;1-P(t)\ge C'\exp(-\la' t^4),\;\;t\ge 0;\quad C',\la'>0;
\quad P(t)=\int_{-\infty}^t p(t')\,dt'.\cr}
$$}

Observe that Theorem 1.1 is a particular case of Theorem 1.2,
so we need to prove only Theorem 1.2.
The plan of the remainder of the paper is as follows.
In Section 2 we present a theorem of Colin de Verdi\`ere and show
how with the help of this theorem to reduce Theorem 1.2
to a lattice--point problem for the Bohr--Sommerfeld
quasi--classical approximation. Section 3 is auxiliary:
here we investigate the asymptotics at infinity of the
Fourier transform of the characteristic function of
a plane domain with non--degenerate inflection points and angular points.
In Section 4 we prove the lattice--point version of
Theorem 1.2 for the number of lattice points inside
a dilated plane oval with finitely many points of
inflection. Section 5 is technical, and here we prove
some lemmas used in Section 4. In Section
6 we prove Theorem 1.2.
Finally in an Appendix
we prove a 2/3--estimate for ovals with semicubic
singularity. This estimate is used in the main part of the paper.

\beginsection 2. Quasi--Classical Approximation \par

Colin de Verdi\`ere proved in [CdV2] the following result:

{\bf Theorem CdV.} {\it If $M$ is a surface of revolution
of simple structure then
$$
\text{\rm Spectrum}\,(-\De)=\{
E_{kl}=Z(k+(1/2),l);\;\; k,l\in\Z,\; |l|\le k\}
\eqno (2.1)
$$
with
$Z(p)=Z(p_1,p_2)\in C^\infty(\R^2)$ such that
$$
Z(p)=Z_2(p)+Z_0(p)+O(|p|^{-1}),\quad
|p|\to\infty,\eqno (2.2)
$$
where
$$
Z_2(p),Z_0(p)\in C^\infty(\R^2\setminus\{ 0\});\quad
Z_2(p)>0,\quad p\not=0,
$$
and
$$
Z_j(\la p)=\la^jZ_j(p),\quad
\forall\,\la>0,p\in\R^2,\quad j=0,2.
$$
In addition, in the sector $\{p_1\ge|p_2|\}$, $Z_2(p)$
satisfies the equation
$$
\PI\int_a^b\sqrt{Z_2(p)-p_2^2f^{-2}(s)}ds=p_1-|p_2|,
\eqno (2.4)
$$
where $a,b$ are the turning points, i.e.,
$$
Z_2(p)-p_2^2f^{-2}(s)=0
\quad\text{\rm for}\quad s=a,b.\eqno (2.5)
$$}

It is to be noted that the Bohr--Sommerfeld
quantization rule is
$$
\PI\int_a^b\sqrt{E_{kl}-l^2f^{-2}(s)}ds
=k+(1/2)-|l|,
\quad k\ge |l|,\eqno (2.6)
$$
which is equivalent to the approximation
$$
E_{kl}=Z_2(k+(1/2),l).\eqno (2.7)
$$
(2.1) implies that
$$
N(R)=\#\{ E_{kl}\le R^2\}
=\#\{(k,l)\: Z(k+(1/2),l)\le R^2,\;|l|\le k\}.
\eqno (2.8)
$$
Define
$$
N_{\BS}(R)=
\#\{(k,l)\: Z_2(k+(1/2),l)\le R^2,\;|l|\le k\}
\eqno (2.9)
$$
(BS stands for Bohr--Sommerfeld).

{\bf Theorem 2.1.}
$$
\liT\IO|N(R)-N_{\BS}(R)|^2R^{-1}dR=0.
\eqno (2.10)
$$

{\it Proof.} Define
$$
M=\{ n=(k+(1/2),l)\: k,l\in\Z;\;|l|\le k\}.
\eqno (2.11)
$$
Then
$$
N(R)-N_{\BS}(R)=
\sum_{n\in M}(\chi(n;R)-\chi_{\BS}(n;R)),
\eqno (2.12)
$$
where $\chi(p;R)$ and $\chi_{\BS}(p;R)$ are the characteristic functions
of the domains
$\{Z(p)\le R^2\}$
and $\{Z_2(p)\le R^2\}$, respectively. Hence
$$
\IO|N(R)-N_{\BS}(R)|^2R^{-1}dR=\sum_{n,n'\in M}I(n,n')
\eqno (2.13)
$$
with
$$
I(n,n')=\IO
(\chi(n;R)-\chi_{\BS}(n;R))
(\chi(n';R)-\chi_{\BS}(n';R))R^{-1}dR.\eqno (2.14)
$$
Now, $\ex C>0$ such that
$$
\chi(n;R)-\chi_{\BS}(n;R)=0\quad\text{if}\quad
\dist\{n,\G_R\}>C|n|^{-1},
\eqno (2.15)
$$
where
$$
\G_R=\{ p\in \R^2\: Z_2(p)=R^2\}=R\G,\quad\quad
\G=\{p\in \R^2\: Z_2(p)=1\}.\eqno (2.16)
$$
Indeed,
$$
\sup_{p\in\R^2}|Z(p)-Z_2(p)|=C_0<\infty,\eqno (2.17)
$$
so if $Z(p)=R^2$ then for $|\ep|<1/2,$
$$
|Z_2(p(1+\ep))-R^2|
\ge |Z_2(p(1+\ep))-Z_2(p)|-|Z_2(p)-Z(p)|
\ge C_1|p|^2\ep-C_0.
$$
Hence if $\ep=C|p|^{-2}$ with
$C=2C_0C_1^{-1}$, then
$$
\pm(Z_2(p(1\pm \ep))-R^2)\ge CC_1-C_0=C_0>0,
$$
and $\ex |\ep_0|<\ep$ such that
$Z_2(p(1+\ep_0))=R^2,$ $|p\ep_0|\le C|p|^{-1}$,
which implies (2.15).

Define $X(p)=(Z_2(p))^{1/2}$. Then $X(p)>0$
and
$$
X(\la p)=\la X(p)\quad\forall\;\la>0,p\in\R^2.
$$
(2.15) implies that
$$
I(n,n')=0\quad\text{if}\quad |X(n)-X(n')|>C_0|n|^{-1}.
\eqno (2.18)
$$
In addition, (2.14), (2.15) imply that
$$
I(n,n')=0\quad\text{if}\quad X(n)\ge 2T,
\eqno (2.19)
$$
and $\forall\; n,n'$,
$$
|I(n,n')|\le CT^{-1}|n|^{-2}.\eqno (2.20)
$$
{}From (2.18)--(2.20)
$$
\sum_{n,n'\in M}I(n,n')\le C\sum_{n\in M\: X(n)\le 2T}
T^{-1}|n|^{-2}\sum_{n'\in M\:|X(n)-X(n')|\le C_0|n|^{-1}}1.
$$
Due to the 2/3--estimate of Sierpinski--Landau--Randol--Colin de Verdi\`ere
[Sie],
[Lan], [Ran], [CdV1],
$$
\sum_{n'\in M\:|X(n)-X(n')|\le C_0|n|^{-1}}1\le C|n|^{2/3},
$$
hence
$$
\sum_{n,n'\in M}I(n,n')\le
CT^{-1}\sum_{n\in M\: X(n)\le 2T}|n|^{-4/3}
\le C_0 T^{-1/3}.\eqno (2.21)
$$
Since (2.13) and (2.21) imply (2.10), Theorem 2.1 is proved.

By (2.9) $N_{\BS}(R)$ is the number of lattice points $(k+(1/2),l)$
of a shifted square lattice
in the sectorial domain
$$
\Om(R)=\{ p\in\R^2\: Z_2(p)\le R^2,\; |p_2|\le p_1\}.
$$
Observe that $\Om(R)$ is $\Om=\Om(1)$ dilated
with the coefficient $R$, so the problem
of finding the asymptotics of $N_{\BS}(R)$
as $R\to\infty$
reduces to a lattice--point problem
about the asymptotics of the number
of lattice points inside $R\Om$,
where $\Om$ is a sectorial domain between
diagonals $p_2=\pm p_1$ which is bounded
by the curve $\G=\{ Z_2(p)=1\}$. By (2.4)
$\G$ is the graph of the function
$$
p_1=|p_2|+\PI\int_a^b\sqrt{1-p_2^2f^{-2}(s)}ds.
$$
Fig.3 shows the curve $\G$
found with the help of computer for
the surfaces of revolution
presented on Fig.2 above.

Sections 3--5 below are devoted to the lattice--point problem
for arbitrary sectorial domain between diagonals which is
bounded by a ``generic'' curve $\G$.

\beginsection 3. Asymptotics of Fourier Transform
of Nonconvex Domains \par

In this section we consider the following
auxiliary problem.
Let $\Om$ be a sectorial domain on a plane, which is
bounded by two segments $[0,z_0],\,[0,z_1]$ and by a
smooth curve $\G$ which goes from $z_0$ to $z_1$ (see
Fig.4). We will assume that $\G$ has at most finitely
many points of inflection where the curvature $\sg(p)$
vanishes, and all the points of inflection are
non--degenerate, i.e., $\di{d\sg\over ds}\not=0$
at these points. We are interested in the asymptotics of the
Fourier transform
$$
\tc(\xi)=\int_{\Om}\exp(ip\xi)\,dp,
$$
as $|\xi|\to\infty$.
With the help of a partition of unity this problem
can be reduced to a series of local problems of the
following type.

Let $\G$ be a smooth curve near some point $p^0\in\G$
and $\chi(p)=1$ from one side of $\G$ and $\chi(p)=0$
from the other side of $\G$ (locally). Let $\f(p)
\in C_0^\infty$ be a $C^\infty$--function with compact support
near $p^0$. Assume in addition that $\f(p)=1$ in
the vicinity of $p^0$.
Then the local problem is: what is the asymptotics of
$$
\tc(\xi)=\int_{\R^2}\f(p)\chi(p)\exp(ip\xi)\,dp,
$$
as $|\xi|\to\infty$?

We will be also interested in the case when $\G$ has an angular
point at $p^0$, and we will consider in a sequence the
following cases: (a) $\G$ is either convex or concave near $p^0$;
(b) $\G$ has an inflection point at $p^0$;
(c) $\G$ has an angular point at $p^0$ with one side which is
either convex or concave and with the other side which is a
straight ray; (d) $\G$ is an angle between two straight rays;
and (e) $\G$ is a straight line (see Fig.5).

In the case (a) the answer is the following well--known
lemma. Let $\G_0\subset\G$ be an open arc on $\G$ such that
$$
p^0\in\G_0\subset\{ p\in\G\: \f(p)=1\},
$$
and
$$
V=\{\xi\in\R^2\setminus\{0\}\:\ex p=p(\xi)\in\G_0\;\text{such that}\;
n_\G(p)=|\xi|^{-1}\xi\},
$$
where $n_\G(p)$ is the vector of normal to $\G$ at $p\in\G$ which
looks in the direction where $\chi(p)=0$. Define
$$
Y(\xi)=\xi\cdot p(\xi),\quad\xi\in\R^2\setminus\{0\}.\eqno (3.1)
$$
For the sake of brevity we will denote $\sg(p(\xi))$ by
$\sg(\xi)$. We will assign a sign to the curvature
$\sg(p)$, so that $\sg(p)>0$ if the region $\{\chi(p)=1\}$
is convex near $p\in\G$, and $\sg(p)<0$ if this
region is concave near $p$. If $p_2=f(p_1)$ is the
equation of $\G$ near some $\hat p\in\G$ in the coordinate
system with an orthonormal basis $e_1,e_2$ such that
$e_2=n_\G(\hat p)$, then
$$
\sg(\hat p)=-f''(\hat p_1).\eqno (3.2)
$$

{\bf Lemma 3.1} (see [Hla]). {\it If $\sg(p)\not=0$ on
$\G_0$ then
$$
\tc(\xi)=(2\pi)^{-1/2}|\xi|^{-3/2}|\sg(\xi)|^{-1/2}
\cos (i(Y(\xi)-\phi))
+O(|\xi|^{-5/2}),\quad \text{\rm if}\;\xi\in V,\eqno (3.3)
$$
where
$$
\phi=(\pi/2)+(\pi/4)\sgn\sg(\xi).\eqno (3.4)
$$}

Assume now that $p^0$ is a non--degenerate point of inflection, i.e.,
$\sg(p^0)=0$ and $\di{d\sg\over ds}(p^0)\not=0$. In this
case $p^0$ is the turning point for $n_\G(p)$, so that for
small deviations of the direction of $\xi$ from $n_\G(p^0)$
we have either two points $p=p_\pm(\xi)$ with
$n_\G(p)=|\xi|^{-1}\xi$ or no such point at all. Define for $\xi
\in V$,
$$
\t(\xi)=\;\;\eqalign{
&1\;\text{if}\;\ex p\in\G_0\;\text{with}\;n_\G(p)=|\xi|^{-1}\xi,\cr
&0\;\text{otherwise.}\cr}\eqno (3.5)
$$
To be definite in the choice of $p_\pm(\xi)$ we will assume that
$\sg(p_+(\xi))>0$ and $\sg(p_-(\xi))<0$. Let
$\sg_\pm(\xi)=\sg(p_\pm(\xi))$ and $Y_\pm(\xi)=\xi\cdot p_\pm(\xi)$.
Let $\a(\xi)$ be the angular coordinate of $\xi$. Let finally
$$
\Ai(t)=\int_{-\infty}^\infty\exp(it\tau+i\tau^3/3)d\tau
$$
be the Airy function. Recall that $\Ai (t)\in C^\infty(\R^1)$ and
when $t\to\infty$,
$$\eqalign{
&\Ai(-t)=\pi^{-1/2}t^{-1/4}\cos(\zeta-(\pi/4))+O(t^{-7/4}),\cr
&\Ai'(-t)=\pi^{-1/2}t^{1/4}\sin(\zeta-(\pi/4))+O(t^{-5/4}),
\quad\zeta=(2/3)t^{3/2},\cr
&|\Ai(t)|,\;|\Ai'(t)|\le C\exp(-t),\cr}
$$
which implies that for $t\not=0$,
$$\eqalign{
&|\Ai(t)-\t(-t)\pi^{-1/2}|t|^{-1/4}\cos(\zeta-(\pi/4))|
\le C|t|^{-7/4},\cr
&|\Ai'(t)-\t(-t)\pi^{-1/2}|t|^{1/4}\sin(\zeta-(\pi/4))|
\le C|t|^{-5/4},\cr}\eqno (3.6)
$$
where $\t(t)=1$ for $t\ge 0$ and $\t(t)=0$ for
$t<0$.

{\bf Lemma 3.2.} {\it Assume that $\sg(p^0)=0$
and $\di{d\sg\over ds}(p^0)\not=0$. Then there
exist real valued functions $a(\a),\,b(\a)$ near
$\a_0=\a(n_\G(p^0))$ such that $a(\a_0)=0,\;
a'(\a_0)\not=0,\;b(\a_0)=Y(n_\G(p^0))$ and
$$\eqalign{
\tc(\xi)&=-i\exp\Bigl(i|\xi|b(\a(\xi))\Bigr)
\Bigl\{|\xi|^{-4/3}\Ai\bigl(|\xi|^{2/3}a(\a(\xi))\bigr)u(\a(\xi))\cr
&+|\xi|^{-5/3}\Ai'\bigl(|\xi|^{2/3}a(\a(\xi))\bigr)v(\a(\xi))\Bigr\}
+O(|\xi|^{-7/3}),\;\;\text{if}\;\xi\in V,\cr}
\eqno (3.7)
$$
where $u(\a),\,v(\a)$ are $C^\infty$ functions and
$$
u(\a_0)=\left|{1\over2}\,{d\sg\over ds}(p^0)\right|^{1/3}.
\eqno (3.8)
$$}

{\bf Corollary.} {\it If $\a(\xi)\not= \a_0$ then
$$
|\chi(\xi) -\t(\xi)(\chi^+(\xi)+\chi^-(\xi))|\le
C|\xi|^{-5/2}|\a(\xi )-\a_0|^{-7/4} \eqno(3.9)
$$
where
$$
\chi^\pm(\xi)=|\xi |^{-3/2}|2\pi \sg _\pm
(\xi)|^{-1/2}\exp(i(Y_\pm (\xi)-\phi_\pm )), \eqno(3.10)
$$
$\phi_\pm=(\pi/2)\pm(\pi/4) $.

 Proof of Corollary}. From (3.6) and (3.7) we obtain that for
$\a(\xi )=\a\not=\a_0$,
$$\eqalign {
\chi _m(\xi)&=-i\exp\bigl(i|\xi |b(\a)\bigr)
|\xi|^{-3/2}\t(\xi)
\Bigl[u_0(\a)\cos\bigl(|\xi|a_0(\a)-(\pi/4)\bigr)+v_0(\a)
\sin \bigl(|\xi |a_0(\a)-(\pi/4)\bigr)\Bigr]\cr
&+O(|\xi|^{-5/2}|\a-\a_0|^{-7/4})\cr} \eqno (3.11)
$$
with some $a_0(\a),u_0(\a)$ and $v_0(\a)$.
On the other hand, Lemma 3.1 gives us that for a fixed
$\a(\xi)=\a\not=\a_0$,
$$
\chi_m(\xi)=\t(\xi)(\chi_m^+(\xi)+\chi_m^-(\xi)) + O(|\xi|^{-5/2}),
\quad |\xi|\to\infty \eqno(3.12)
$$
Comparing (3.11) with (3.12) we obtain that the main terms
in these two asymptotics coincide, hence
$$
\chi_m(\xi)=\t(\xi)(\chi_m^+(\xi)+\chi_m^-(\xi)) + O(|\xi|^{-5/2}|\a-
\a_0|^{-7/4}),
$$
which proves Corollary.

{\it Proof of Lemma 3.2.} Consider an orthonormal basis $e_1,e_2$ on the
plane with $e_2=n_\G (p^0)$. Let
$
p=p^0+p_1e_1+p_2e_2\,$, $\xi=\xi_1e_1+\xi_2e_2\,$ and $p_2=f(p_1)$
be the equation of $\G$ near $p^0$. Observe that $f(0)=f'(0)=f''(0)=0$
and we can choose the direction of $e_1$ in such a way that
$f'''(0)=\di{d\sg\over {ds}}(p^0).$ Integrating by parts in $p_2$ we obtain
that
$$\eqalign {
    \tilde \chi (\xi )
&=\int_{\R^2} \f(p)\chi(p)\exp(ip\xi) dp=
-i\xi_2^{-1}\exp(ip^0\xi )\int_{-\infty}^{\infty }\f(p_1,f(p_1))
\exp(ip_1\xi_1+if(p_1)\xi_2)dp_1 \cr
&+O(|\xi|^{-N})
=-i\xi_2^{-1}\exp(ip^0\xi )\int_{-\infty}^{\infty }\f(p_1,f(p_1))
\exp[i|\xi| \Phi(p_1,\a(\xi))]dp_1 +O(|\xi|^{-N}),\cr }
$$
where
$\Phi(p_1,\a)=p_1\sin (\a-\a_0 )+f(p_1)\cos(\a -\a_0) $
is a $C^\infty $ real valued function with
$$\Phi={\partial \Phi \over {\partial p_1}}
={\partial^2 \Phi \over {\partial p_1^2}}=0,\quad
{\partial^3 \Phi \over {\partial p_1^3}}={d\sg\over {ds}}(p^0),\quad
\text { and } {\partial^2\Phi\over{\partial p_1\partial \a}}\not= 0,$$
at $p_1=0,\a=\a_0.$ (3.7) follows now from Theorem 7.7.18 in [H\"or2].
Lemma 3.2 is proved.

Let us turn now to angular points. So, assume that $\G$ is a smooth
curve near $p^0\in \G$ and $L$ is a straight line which intersects
$\G$ at $p^0$ transversally.
Then near $p^0$, $\G$ and $L$ divide the plane into four parts. Let
$\chi (p)=1$ in one of these parts and $\chi (p)=0$ in the remainder.
Again we are interested in asymptotics of $\tilde \chi (\xi)$ as
$|\xi |\to\infty.$

Define the auxiliary functions
$$
P_\pm (y)=(2\pi)^{-1/2}\exp(\pm\pi i/4)\int_{-\infty}^y
\exp(\mp  it^2/2)dt, \eqno(3.13)
$$
which are $C^\infty $ bounded functions on $\R^1$ such that
$$ \eqalign{
&P_\pm(-\infty )=0,\quad P_\pm (0)=1/2,\quad P_\pm(\infty)=1;\cr
&P_+(y)=\ar {P_-(y)},\quad P_\pm (y)+P_\pm (-y)=1;\cr
&|P_\pm (y)-\t (y)|\le C (1+|y|)^{-1}\cr } \eqno (3.14)
$$

{\bf Lemma 3.3.} {\it Assume $\sg (p)\not= 0$ in vicinity of $p^0\in \G$
and $L$ intersects $\G$ at $p^0$ transversally. Then there exists
$C^\infty $ real valued functions $a(\a )$ and $b(\a )$
 near $\a_0=\a (n_\G(p^0))$
such that
$a(\a_0)=0, a'(\a_0)\not=0, b(\a_0)=Y(n_\G(p^0))$
and
$$
\tilde \chi (\xi)=|\xi|^{-3/2}\exp[i|\xi|b(\a(\xi))-i\phi]
P_\zeta\bigl(|\xi|^{1/2}a(\a(\xi))\bigr)u(\a(\xi))+O(|\xi|^{-2}),
\;\text { if } \xi\in V,\eqno(3.15)
$$
where $\phi =(\pi/2)+(\pi/4)\sgn \sg(p^0)$,
$\zeta=\sgn\sg(p^0)$ and
$u(\a)$ is a $C^\infty $ function near $\a_0$ with
$u(\a_0)=|2\pi\sg(p^0)|^{-1/2}$.}

{\bf Corollary.} {\it If $\xi\in V$ and $\a(\xi)\not=\a_0=\a(n_\G(p^0))$,
then
$$
\Bigl|\tilde \chi(\xi)- |\xi|^{-3/2}|2\pi\sg(\xi )|^{-1/2}\t(\xi)
\exp[i(Y(\xi)-\phi)]\Bigr|\le C|\xi|^{-2}|\a(\xi)-\a_0|^{-1},\eqno(3.16)
$$
where $\t(\xi)=1 $ if $\ex p=p(\xi )\in \G$
with $n_\G(p)=|\xi|^{-1}|\xi|$, and $\t (\xi )=0$
otherwise.

Proof of Corollary.} From (3.14), (3.15), we obtain that for
$\a(\xi)=\a\not =\a_0$
$$
\tilde \chi (\xi)=|\xi|^{-3/2}\exp [i|\xi| b(\a)-i\phi ]\t(\xi )u(\a)
+O(|\xi|^{-2}|\a-\a_0|^{-1}), \quad\xi\in V
$$
On the other hand, when $\a(\xi )=\a\not =\a_0$ is fixed, Lemma 3.1
proves that
$$
\tilde \chi(\xi)=|\xi|^{-3/2}|2\pi \sg(\xi)|^{-1}\t(\xi)\exp[i(Y(\xi)-
\phi )] +O(|\xi |^{-2}), \,\; |\xi|\to\infty,\;\; \xi\in V.
$$
Comparing these two asymptotics we conclude that the main terms in them
coincide, and therefore (3.16) holds.

{\it Proof of Lemma 3.3.} Consider a basis $e_1,e_2$ on the plane,
where $e_1$ is a unit tangent vector to $\G$
at $p^0$ and $e_2$ is a unit tangent vector
to $L$ at $p^0$.
Let $e_1^\perp,e_2^\perp$ be a dual basis,
$e_i^\perp \cdot e_j=\de_{ij}$,
and $p=p^0+p_1e_1+p_2e_2,\, \xi=\xi_1e^\perp+\xi_2e^\perp_2$.
Then
$$\eqalign {
\tilde \chi(\xi)&=\int_{\R^2}\f(p)\chi(p)\exp(ip\xi)dp\cr
&=J\exp(ip^0\xi)i\xi_2^{-1}\int_0^\infty\f(p_1,f(p_1))
\exp(ip_1\xi_1+if(p_1)\xi_2)dp_1+O(|\xi|^{-N})\cr}\eqno(3.17)
$$
where $p_2=f(p_1)$ is the equation of $\G$ and $J$ is the Jacobian.
Observe that
$ p_1\xi_1+f(p_1)\xi_2=|\xi|\Phi(p_1;\a(\xi)) $
where $\Phi(p_1,\a)$ is a $C^\infty$ real valued function with
$$
\Phi={\partial\Phi\over {\partial p_1}}=0 ,\quad
{\partial^2\Phi\over {\partial p_1^2}}\not= 0 ,\quad
{\partial^2\Phi\over {\partial p_1\partial \a}}\not= 0,$$
at $p_1=0,\a=\a_0$.

By theorem 7.5.13 in [H\"or2] there exists a $C^\infty $
change of variable $t=t(p_1,\a)$ with $t(0,\a_0)=0$,
$\di{\partial t\over {\partial p_1 }}(0,\a_0)\not= 0$
such that
$\Phi(p_1,\a)=-\zeta(t^2/2)+b(\a)$, $b(\a_0)=0$,
where $\zeta =\sgn\sg(p^0)$. Hence
$$\eqalign {
I(\xi)&\equiv \int_0^\infty \f(p_1,f(p_1))\exp(i|\xi|\Phi(p_1,\a))
dp_1\cr
&= \exp(i|\xi|b(\a))\int_{t(0,\a)}^\infty
\psi (t;\a)\exp(-i|\xi|\zeta t^2/2)
{dp_1\over {dt}}dt,\cr }
$$
where $\psi(t,\a)=\f(p_1,f(p_1)),\;p_1=p_1(t,\a)$.
Observe that $ \psi(t,\a)=1 $
near $(0,\a_0)$.  Let
$$
{dx_1\over {dt}}=u_0(\a)+tu_1(\a,t)
$$
Then
$$
I(\xi)=\exp(i|\xi|b(\a))u_0(\a)|\xi|^{-1/2}
P_\zeta(-t(0,\a)|\xi|^{1/2})+O(|\xi|^{-1}) \eqno (3.18)
$$
{}From (3.17), (3.18) we obtain (3.15). Lemma 3.3 is proved.

{\bf Lemma 3.4.} {\it Let $\chi(p)=1$ between two rays
$L_1=\{p^0+\la e_1,\;\la\ge 0\}$
and $L_2=\{p^0+\la e_2,\;\la\ge 0\}$,
where $e_1,e_2$ are linearly independent vectors, and $\chi(p)=0$
otherwise. Then
$$
|\tc(\xi)|\le C|\xi|^{-1}\bigl[(1+|\xi\cdot e_1|)^{-1}+(1+|\xi\cdot
e_2|)^{-1}\bigr],\eqno (3.20)
$$
and for $\la\to\pm\infty$,
$$
\tc(\la e_j^\perp)=C_j\la^{-1}+O(|\la|^{-2}),\quad j=1,2;
\quad e_i^\perp\cdot e_j=\de_{ij}.
\eqno (3.21)
$$

Proof.} Let $p=p^0+p_1e_1+p_2e_2,\;\xi=\xi_1 e_1^\perp
+\xi_2 e_2^\perp$. Without loss of generality we may
assume that $|\xi_2|\ge|\xi_1|$. Then
integrating by parts in $p_2$ we obtain that
$$
\tc(\xi)=i\xi_2^{-1}\exp(ip^0\xi)
\int_0^\infty\f(p_1,0)\exp(ip_1\xi_1)dp_1+O(|\xi|^{-2}).
$$
Integrating now by parts in $p_1$ we obtain (3.20), while
setting $\xi_1=0$ we obtain (3.21). Lemma 3.4 is proved.

{\bf Lemma 3.5.} {\it Let $\chi(p)=1$ from one side of the line
$L=\{ p^0+\la e,\;\la\in\R\}$ and $\chi(p)=0$ from the other side
of $L$. Then
$$
|\tc(\xi)|\le C|\xi|^{-1}(1+|\xi\cdot e|)^{-1},\quad
\tc(\la e^\perp)=C\la^{-1}\exp(i\la p^0\cdot e^\perp)+
O(|\la|^{-2}),\quad |\la|\to\infty.
$$}

The proof of Lemma 3.5 is similar to the proof of Lemma 3.4.

\beginsection 4. Lattice--Point Problem \par

Let $Z(p)\in C^\infty(\R^2\setminus\{ 0\})$ be
a $C^\infty$ positive function homogeneous of degree 2. Define
$$\eqalign{
&N(R)=\#\{ n=(n_1,n_2)\in\Z^2\:
Z(n_1+(1/2),n_2)\le R^2,\;|n_2|\le n_1\},\cr
&\Om(R)=\{ p\in\R^2\: Z(p)\le R^2,\; |p_2|\le p_1\},\cr
&A(R)=AR^2=\Area\Om(R).\cr}\eqno (4.1)
$$
The lattice--point problem we are interested
in is to evaluate $N(R)-A(R)$ as $R\to\infty.$

Let
$$
\G=\{ p\in\R^2\: Z(p)=1,\;|p_2|\le p_1\},
$$
and $z_0,z_1\in\G$ be the endpoints of $\G$ with
$z_{01}=-z_{02}>0$ and $z_{11}=z_{12}>0$. For
$p\in\G$ denote by $n_\G(p)$ the vector of outer
normal to $\G$ at $p$. Observe that
$$
n_\G(p)=|\grad Z(p)|^{-1}\grad Z(p),
$$
and
$$
p\cdot\grad Z(p)=2Z(p)>0,
\eqno (4.2)
$$
hence
$$
p\cdot n_\G(p)=2|\grad Z(p)|^{-1}\grad Z(p)>0.
\eqno (4.3)
$$
Denote by $\sg(p)$ the curvature of $\G$
at $p\in\G$ with a sign, so that $\sg(p)>0$
if $\Om$ is convex near $p$, $\sg(p)<0$
if $\Om$ is concave near $p$.
In what follows we assume the following

{\bf Hypothesis D.} {\it (i) $\sg(p)\not=0$
everywhere on $\G$ except, maybe, a finite set
$W=\{w_1,\dots, w_K\},$ $z_0,z_1\not\in W$, and
$$
{d\sg\over ds}(w_k)\not=0,\quad k=1,\dots, K,
\eqno (4.5)
$$
where $s$ is the natural coordinate on $\G$.
(ii) $\forall\, w_k\in W$ the vector $\nu_k=n_\G(w_k)$
is either rational, i.e., $n\cdot\nu_k=0$ for some
$n\in\Z^2,\;n\not=0$, or Diophantine in the sense that
$\ex 1>\zeta>0$ and $C>0$ such that
$$
|n\cdot\nu_k|>{C\over|n|^{1+\zeta}},\quad
\forall\, n\in\Z^2,\;n\not=0.
\eqno (4.6)
$$}

Without loss of generality we
may assume that
$$
0=s(z_0)<s(w_1)<\dots<s(w_K)<s(z_{1}),
\eqno (4.7)
$$
where $s(p)$ is the natural coordinate of $p\in\G$.

We call $\xi\in\R^2\setminus\{0\}$ rational if
$n\cdot\xi=0$ for some $n\in\Z^2\setminus\{0\}$.
It is to be noted that $\xi$ is rational
iff the set
$$
L(\xi)
=(\Z^2\setminus\{0\})\cap \{\la\xi,\,\la\in\R\}
$$
is non--empty. Let
$$
\G_{\rat}=\{ p\in\G\: n_\G(p)\quad\text{is rational}\}.
$$
For $p\in\G$ define
$$
Y(p)=p\cdot n_\G(p).\eqno (4.8)
$$
By (4.3) $Y(p)>0$.

{\bf Theorem 4.1.} {\it Assume Hypothesis D holds. Then
$$
N(R)=AR^2+R^{2/3}\sum_{k\:\nu_k\;\text{\rm is rational}}
\Phi_k(R)+R^{1/2}F(R),\eqno (4.9)
$$
where $\Phi_k(R)$ are continuous periodic functions,
$$
\Phi_k(R)=(1/2)3^{-2/3}\G(2/3)\pi^{-4/3}\left|
{d\sg\over ds}(w_k)\right|^{-1/3}
\sum_{n\in L(\nu_k)} (-1)^{n_1}|n|^{-4/3}\sin(2\pi Y(w_k)|n|R)
\eqno (4.10)
$$
and $F(R)\in B^2$. The Fourier series of $F(R)$ is
$$
F(R)=\PI\sum_{p\in\G_{\rat}\setminus W}\t(p)|\sg(p)|^{-1/2}
\sum_{n\in L(n_\G(p))}(-1)^{n_1}|n|^{-3/2}
\cos (2\pi Y(p)|n|R-\phi(p)),\eqno (4.11)
$$
where
$$
\t(p)=\;\;
\eqalign{
&1\quad\text{\rm if}\quad p\not= z_0,z_1,\cr
&(1/2)\quad\text{\rm if}\quad p=z_0,z_1,\cr}
$$
and $\phi(p)=(\pi/2)+(\pi/4)\sgn\sg(p)$.}

We need Theorem 4.1 to prove our main Theorem 1.2. For the
sake of completeness we want to formulate another theorem
which is not needed for the proof of Theorem 1.2,
but which is of interest by itself.

Let  $\a\in\R^2$ be a fixed point on the plane. Define
$$
N(R;\a)=\#\{ n\in\Z^2\: Z(n+\a)\le R^2\}.
$$

{\bf Theorem 4.2.} {\it Assume Hypothesis D holds. Then
$$
N(R;\a)=AR^2+R^{2/3}\sum_{k\:\nu_k\;\text{\rm is rational}}
\Phi_k(R;\a)+R^{1/2}F(R;\a),
$$
where $\Phi_k(R;\a)$ are continuous periodic functions of $R$,
$$
\Phi_k(R;\a)=(1/2)3^{-2/3}\G(2/3)\pi^{-4/3}\left|
{d\sg\over ds}(w_k)\right|^{-1/3}
\sum_{n\in L(\nu_k)} |n|^{-4/3}\sin(2\pi (Y(w_k)|n|R-n\cdot\a))
$$
and $F(R;\a)\in B^2$ in $R$.
The Fourier series of $F(R;\a)$ is
$$
F(R;\a)=\PI\sum_{p\in\G_{\rat}\setminus W}|\sg(p)|^{-1/2}
\sum_{n\in L(n_\G(p))}|n|^{-3/2}
\cos (2\pi Y(p)|n|R-\phi(p,n,\a)),\eqno (4.12)
$$
where
$\phi(p,n,\a)=(\pi/2)+(\pi/4)\sgn\sg(p)+2\pi n\cdot\a$.}

Theorem 4.2 is a generalization of Theorem 1.1
in [Ble1] to non--convex  domains.

{\it Proof of Theorem 4.1}. $N(R)$ can be written as
$$
N(R)=\sZ\chi(n+\a;R),\quad\a=(1/2,0),\eqno (4.13)
$$
where $\chi(p;R)$ is the characteristic function of $\Om(R)$.
Define for $\de>0$,
$$
N_\de(R)=\sZ\chi_\de(n+\a;R),\eqno (4.14)
$$
where
$$
\chi_\de(p;R)=\int_{\Om(R)}\de^{-2}\f(\de^{-1}(p-p'))\,dp'
=\chi(\cdot;R)\ast(\de^{-2}\f(\de^{-1}\cdot))(p),
\eqno (4.15)
$$
and
$$
\f(p)\in C_0^\infty(\R^2);\;\; \f(p)=\f_0(|p|);\;\;
\f(p)\ge 0;\;\; \int_{\R^2}\f(p)dp=1;\;\;
\f(p)=0\quad\text{\rm when}\quad |p|\ge 1.
$$

{\bf Lemma 4.3.} For $T>1$,
$$
\IT|N_\de(R)-N(R)|^2R^{-1}dR\le C\de T(T^{-1/3}+\de).
\eqno (4.16)
$$

Proof of this and all subsequent lemmas is given in
the next section. For what follows we put
$$
\de=T^{-1}.\eqno (4.17)
$$
In this case (4.16) reduces to
$$
\IT|N_\de(R)-N(R)|^2R^{-1}dR\le CT^{-1/3}.\eqno (4.18)
$$

By the Poisson summation formula
$$
N_\de(R)-AR^2=R^2\sz\tpd\tc(2\pi nR)(-1)^{n_1},
\eqno (4.19)
$$
where
$$
\tc(\xi)=\int_{\Om}\exp(ip\xi)dp,\quad
\Om=\{ p\in\R^2\: Z(p)\le 1,\;|p_2|\le p_1\}.
$$
Let us consider a partition of unity on the projective line
$\R\bold P^1$,
$$
\sum_{l=1}^L\psi_l(\xi)=1,\quad 0\le\psi_l(\xi)\le 1,\quad\psi_l(\xi)
\in C^\infty(\R\bold P^1),
$$
and lift it to $\R^2\setminus\{0\}$ putting
$$
\psi(\la\xi)=\psi(\xi)\quad\forall\,\la\in\R^1\setminus\{0\}.
$$
Let
$$
\G^{(l)}=\{ p\in\G\: n_\G(p)\in\supp\psi_l(\xi)\},
$$
and $\G_j^{(l)}$, $j=1,\dots,J_l$, be connected components
of $\G^{(l)}$. Without loss of generality we may assume that
each $\G_j^{(l)}$ contains
at most one singular point $z_0,z_1,w_1,\dots,w_K$.

Then, for a given $l$ consider a partition of unity in the
$p$--plane,
$$
\sum_{m=1}^M \chi_m(p)=1,\quad 0\le \chi_m(p)\le 1,\quad
\chi_m(p)\in C^\infty(\Om),
$$
such that for each $\Gjl$, $1\le j\le J_l$, $\ex$
a unique $\chi_m(p)$ which is $\not\equiv 0$ on $\Gjl$.
Without loss of generality we may also assume that
for each $m$, $\chi_m(p)$ contains at most one
singular point.

Define
$$
E_{lm}(R)=R^2\sz\pn\tpd\tc_m(2\pi nR)(-1)^{n_1},
\eqno (4.20)
$$
where
$$
\tc_m(\xi)=\int_{\Om}\chi_m(p)\exp(ip\xi)dp.
$$
Then by (4.19)
$$
N_\de(R)-AR^2=\sum_{l,m}E_{lm}(R).
$$
Now we evaluate $\Elm$. Let $\G_m=\G\cap\supp\chi_m$.

{\bf Lemma 4.4} ($\chi_m$ out of $\G^{(l)}$). {\it
If $\G_m\cap\G^{(l)}=\emptyset$ then
$$
\sup_{1\le R\le T}|\Elm|\le C\log^2T.
$$}

Next we consider the case when $\Gjl\subset \G_m$
and $\G_m$ contains no singular point.
In this case $\forall\,\xi\in
\supp\psi_l$ $\ex$ a unique $p(\xi)\in\Gjl$ such that
$n_\G(p(\xi))=|\xi|^{-1}\xi$.
Denote by
$$
Y_0(\xi)=\xi\cdot n_\G(p(\xi)),\quad \sg_0(\xi)=\sg(p(\xi)).
$$

{\bf Lemma 4.5} (regular $\chi_m$ on $\G$).
{\it Assume that $z_0,z_1,w_1,\dots,w_K\not\in\G_m$ and
$\Gjl\subset\G_m$. Then $\Elm=R^{1/2}\Flm$,
where
$\Flm\in B^2$ and
$$
\Flm=\PI\sz\pn|n|^{-3/2}|\sg_0(n)|^{-1/2}
\cos (2\pi Y_0(n)R-\phi),\eqno (4.21)
$$
where $\phi=(\pi/2)+(\pi/4)\sgn \sg(p)$, $p\in\G_m$.}

Assume now that some inflection point $w_k$ lies
inside $\Gjl$ and $\Gjl\subset\G_m$. Then for every
$\xi\in\supp\psi_l\setminus\{\la\nu_k,\la\in\R\}$
two possibilities exist: either $\ex$
 two points $p_\pm(\xi)
\in\Gjl$ such that $n_\G(p_\pm(\xi))=|\xi|^{-1}\xi$,
or there is no such point at all. Define the function
$\t(\xi)$ which is equal to 1 in the first case,
and which is equal to 0 in the second case.
For the sake of definiteness we will assume
that $\pm\sg(p_\pm(\xi))>0$.
Denote
by
$$
Y_\pm(\xi)=\xi\cdot n_\G(p_\pm(\xi));
\quad \sg_\pm(\xi)=\sg(p_\pm(\xi)).
$$

{\bf Lemma 4.6.} ($\chi_m$ on inflection point).
{\it If $\nu_k=n_\G(w_k)\in\Gjl\subset\G_m$
and $\nu_k$ is rational then
$$
\Elm=R^{2/3}\Plm+R^{1/2}\Flm,
\eqno (4.23)
$$
where $\Plm$ is a periodic continuous function,
which is given in (4.10),
and $\Flm\in B^2$.
The Fourier series of $\Flm$ is
$$
\Flm=\PI\sum_\pm\sum_{n\in\Z^2\setminus\{\la\nu_k,
\,\la\in\R\}}\pn\t(n)
|n|^{-3/2}|\sg_\pm(n)|^{-1/2}\cos(2\pi Y_\pm(n)R-
\phi_\pm),\eqno (4.25)
$$
where $\phi_\pm=(\pi/2)\pm(\pi/4)$.

If $\nu_k$ is Diophantine then $\Elm=R^{1/2}\Flm$
where $\Flm\in B^2$,
and the Fourier series of $\Flm$ coincides with (4.25).}

{\bf Lemma 4.7} (angular $\chi_m$). {\it
If $n_\G(z_0)\in\Gjl\subset\G_m$ then $\Elm=R^{1/2}\Flm$,
where $\Flm\in B^2$, and the Fourier series of $\Flm$ is
$$
\Flm=\PI\sz\pn\t(n)|n|^{-3/2}
|\sg_0(n)|^{-1/2}\cos(2\pi Y_0(n)R-\phi),\eqno (4.26)
$$
where
$$
\t(n)=\;\;\eqalign{
&1\;\text{\rm if}\; n=n_\G(p)\;\text{\rm for some}
\; p\in\Gjl\setminus\{ z_0\},\cr
&(1/2)\;\text{\rm if}\; n=n_\G(z_0),\cr
&0\;\text{\rm otherwise},}\eqno (4.27)
$$
and $\phi=(\pi/2)+(\pi/4)\sgn(z_0)$.}

Proof of Lemmas 4.3--4.7 is given in the next section.

{\it End of the proof of Theorem 4.1.} Let us fix $l$. By Lemma 4.4,
$$
\sum_{m\: \G_m\cap\G^{(l)}=\emptyset}\liT\IT
|R^{-1/2}\Elm|^2dR=0,
$$
and by Lemmas 4.5--4.7,
$$
\sum_{m\: \G_m\cap\G^{(l)}\not=\emptyset}\Elm
=R^{2/3}\sum_{m\in I_1(l)}\Plm
+R^{1/2}\sum_{m\in I_2(l)}\Flm,
$$
where
$$
I_1(l)=\{ m\:\ex k,j\;\text{such that}\;
w_k\in\Gjl\subset\G_m\;\text{and}\;\nu_k
\;\text{is rational}\},
\quad I_2(l)=\{ m\:\G_m\cap\G^{(l)}\not=\emptyset\}.
$$
Making a summation over $l$ we obtain that if we define
$$
\Phi(R)=\sum_{l=1}^L\sum_{m\in I_1(l)}\Plm,\quad
F(R)=\sum_{l=1}^L\sum_{m\in I_2(l)}\Flm,
$$
then
$$
\liT\IT|N_\de(R)-AR^2-R^{2/3}\Phi(R)
-R^{1/2}F(R)|^2R^{-1}dR=0.
$$
Lemma 4.3 implies that the same is true for $N(R)$:
$$
\liT\IT|N(R)-AR^2-R^{2/3}\Phi(R)
-R^{1/2}F(R)|^2R^{-1}dR=0.
$$
Since $\Flm\in B^2$, $F(R)\in B^2$ as well (see [Bes]).
The Fourier series for $\Phi(R)$ and $F(R)$
are the sum of the Fourier series for
$\Plm$ and $\Flm$, respectively. This proves
the formulas (4.11), (4.12). Theorem 4.1 is proved.

Theorem 4.2 is proved in the same way.

\beginsection 5. Proof of Lemmas \par

{\it Proof of Lemma 4.3.} We have:
$$
N_\de(R)-N(R)=\sZ(\chi_\de(n+\a;R)-\chi(n+\a;R)).
$$
The support of
$\chi_\de(p;R)-\chi(p;R)$ is concentrated in the $\de$--neighborhood of
$\partial\Om(R)$.
Observe that $\partial\Om(R)$ consists of the curve $R\G$
and two segments $[0,Rz_0]$ and $[0,Rz_1]$
oriented along the diagonals $p_1\pm p_2=0$. If $\de$ is small, then
the $\de$--neghborhood of $[0,Rz_0]\cup [0,Rz_1]$
does not contain points $n+\a$, $n\in\Z^2$, $\a=(1/2,0)$,
and it does not contribute to $N_\de(R)-N(R)$.

Now,
$$
I\equiv\IO|N_\de(R)-N(R)|^2R^{-1}dR=\sum_{n,n'}I_\de(n,n'),
$$
where
$$
I_\de(n,n')=\IO(\chi_\de(n+\a;R)-\chi(n+\a;R))
(\chi_\de(n'+\a;R)-\chi(n'+\a;R))R^{-1}dR.
$$
Observe that $I_\de(n,n')=0$ unless both $n+\a$ and $n'+\a$
lie in $\de$--neighborhood of $R\G$ for some $R$, hence
$\ex C>0$ such that
$$
I_\de(n,n')=0,\quad\text{if either}\quad
|X(n+\a)-X(n'+\a)|\ge C\de
\quad\text{or}\quad X(n+\a)\ge CT,
$$
where $X(p)=(Z(p))^{1/2}$. In addition, for all $n,n'$,
$$
I_\de(n,n')\le CT^{-1}\de X(n)^{-1}\le C_0T^{-1}\de|n|^{-1}.
$$
Therefore,
$$
I\le\sum_{n\: X(n+\a)\le CT}C_0T^{-1}\de|n|^{-1}
\sum_{n'\:|X(n+\a)-X(n'+\a)|\le C\de}1
$$
By 2/3--estimate (see, e.g., [CdV1])
$$
\sum_{n'\:|X(n+\a)-X(n'+\a)|\le C\de}1
\le C(n^{2/3}+\de|n|),
$$
hence
$$
I\le C_1T^{-1}\de\sum_{n\: X(n+\a)\le CT}
|n|^{-1}(|n|^{2/3}+\de|n|)
\le C_2\de T(T^{-1/3}+\de).
$$
Lemma 4.3 is proved.

{\it Proof of Lemma 4.4.} Let us consider two cases:
 when $0\in\supp \chi_m(p)$ and when $z_0\in\G_m$
(or $z_1\in\G_m$)
and $\G_m\cap\G^{(l)}=\emptyset$;
all the other cases are simpler. Due to (4.20),
$$
E_{lm}=\sz K(n)
$$
with
$$
K(n)=\tp(2\pi n\de)\psi_l(n)R^2\tc_m(2\pi Rn)
(-1)^{n_1}.\eqno (5.1)
$$
By Lemma 3.4, if $0\in\supp\chi_m(p)$ then
$$
R^2|\tc_m(R\xi)|
\le C|\xi|^{-1}\max\{
R^{-1}+|\xi_1+\xi_2|^{-1},
R^{-1}+|\xi_1-\xi_2|^{-1}\},
$$
so if $n\in\La$, where
$$
\La=\{n\in\Z^2\:|n_1+n_2|\ge|n_1-n_2|>0\},
$$
then
$$
R^2|\tc_m(2\pi Rn)|\le C|n|^{-1}|n_1-n_2|^{-1}.
$$
Therefore $ I_0\equiv\sum_{n\in\La} K(n) $
is estimated as
$$
I_0\le C\sum_{n\in\La}|\tpd|\cdot|n|^{-1}
|n_1-n_2|^{-1}.
$$
Since
$\f(p)\in C_0^\infty(\R^2)$,
$$
|\tp(\xi)|\le C(1+|\xi|)^{-5},\eqno (5.2)
$$
and
$$
\sum_{|n|>T^{2}} |\tpd|\le C_0,\quad \de=T^{-1}.
$$
On the other hand
$$
\sum_{n\in\La,\,|n|\le T^2}|n|^{-1}|n_1-n_2|^{-1}
\le C\log^2T,\eqno (5.3)
$$
Hence $I_0\le C\log^2T$. The sum over
$\{|n_1-n_2|\ge|n_1+n_2|>0\}$ is estimated similarly
hence we obtain
$$
\left|\sum_{n\: n_1-n_2\not=0,\,n_1+n_2\not=0}
K(n)\right|\le C\log^2T.
$$
By Lemma 3.4
$$
\sum_{n\: n_1=n_2\not=0}K(n)
=C_1R\sum_{n\: n_1=n_2\not= 0}\tpd\pn n_1^{-1}(-1)^{n_1}
+O(1).
$$
Since
$\tilde\varphi (\xi)$ and $\psi_l(\xi)$ are even functions,
$$
\sum_{n\: n_1=n_2\not= 0}\tpd\pn n_1^{-1}(-1)^{n_1}
=0,
$$
and thus
$$
\sum_{n\: n_1=n_2\not=0}K(n)=O(1).
$$
Similarly,
$$
\sum_{n\: n_1=-n_2\not=0}K(n)=O(1).
$$
As a result,
$$
|E_{lm}(R)|=\left|\sz K(n)\right|\le C\log^2T,
$$
which was stated.

Assume now that $z_0\in\G_m$ and $\G_m\cap\G^{(l)}=\emptyset$.
Let $\Om_0$ be the angle with the vertex at $z_0$
and the sides which go along $[z_0,0]$ and the
tangent vector to $\G$ at $z_0$, so that
$\Om_0$ is a linear approximation to $\Om$ near $z_0$.
Assume for the sake of definiteness that $\sg(z_0)>0$.
Then $\Om_0\supset\Om$ near $z_0$. Let
$$
\tc_m^{(0)}(\xi)=\int_{\Om_0}\chi_m(p)\exp(ip\xi)dp.
$$
Then
$$
\tc_m^{(0)}(\xi)-\tc_m(\xi)=
\int_{\Om_0\setminus\Om}\chi_m(p)\exp(ip\xi)dp.
$$
is estimated as follows. Integrating in the direction orthogonal
to $\xi$ we obtain that
$$
\tc_m^{(0)}(\xi)-\tc_m(\xi)=
\int_{t_0}^{t_1}\ar\chi_m(t)\exp(it|\xi|)dt,
$$
where $\ar\chi_m(t)$ is equal to zero in vicinity of $t_1$ and
$$
\ar\chi_m(t)=a_1(t-t_0)^2+a_2(t-t_0)^3+\dots
$$
in vicinity of $t_0$. This implies
that $|\tc_m^{(0)}(\xi)-\tc_m(\xi)|\le C(1+|\xi|)^{-3}$,
hence
$$
\sz (K(n)-K^{(0)}(n))=O(R^{-1}),
$$
where
$$
K^{(0)}(n)=\tpd\pn R^2\tc_m^{(0)}(2\pi Rn)(-1)^{n_1}.
$$
Now, $\Om_0$ is an angular domain. Using the same arguments
as in the case $0\in\supp\chi_m$, we obtain
that
$$
\sz K^{(0)}(n)=O(\log^2T).
$$
This proves that
$$
\sup_{1\le R\le T}|\Elm|\le C\log^2T.
$$
Lemma 4.4 is proved.

We omit proof of Lemma 4.5 and pass now
to more complicated Lemmas 4.6, 4.7. Lemma
4.5 is proved similarly, with
some simplifications (see also the proof
of Theorem 3.1 in [Ble1]).

{\it Proof of Lemma 4.6.} Assume $\nu_k\in\G^{(l)}_j\subset\G_m$.
Let us split $E_{lm}(R)$ into two parts:
$$
E_{lm}^{(1)}(R)=\sum_{n\in\Z^2\cap L_\g,\,n\not=0} K(n),\quad
E_{lm}^{(2)}(R)=\sum_{n\in\Z^2\setminus L_\g} K(n),\quad
\eqno (5.4)
$$
where $K(n)$ is defined in (5.1) and
$$
L_\g=\{\xi\in\R^2\:\dist(\xi,L)\le \g\},\quad
L=\{\xi=\la\nu_k,\quad \la\in\R^1\}.
$$
with some $\g>0$ which will be chosen later.
Let us evaluate $E_{lm}^{(1)}(R)$.

Assume first that $\nu_k$ is rational.
Then we can choose $\g>0$ such that $L_\g\setminus L$
contains no integer points. In this case
$$
E_{lm}^{(1)}(R)
=\sum_{n\in L(\nu_k)}
\tpd R^2\tc(2\pi nR)(-1)^{n_1}.
$$
By Lemma 3.5, if $\xi\in L$ then
$$
\tc(\xi)=-i\exp(iY(\xi))|\xi|^{-4/3}\Ai(0)|(1/2)\sg'(w_k)|^{-1/3}
+O(|\xi|^{-5/3}),\quad \Ai(0)=3^{-2/3}\G^{-1}(2/3),
$$
hence
$$
E_{lm}^{(1)}(R)
=R^{2/3}\Phi_{lm}(R;\de)+O(R^{1/3}),
$$
with
$$
\Phi_{lm}(R;\de)=
\sum_{n\i L(\nu_k)}
  \tpd(-i)\sin(2\pi Y(n)R)
\times|2\pi n|^{-4/3}\Ai(0)|(1/2)\sg'(w_k)|^{-1/3}(-1)^{n_1}.
$$
Since
$$
\sum_{q=1}^\infty
|1-\tp(2\pi qr\de)|q^{-4/3}\le C\de^{1/3}=CT^{-1/3},
$$
we obtain that
$$
E_{lm}^{(1)}(R)
=R^{2/3}\Phi_{lm}(R)+O(R^{1/3}+R^{2/3}T^{-1/3}),\eqno (5.5)
$$
where
$$
\Phi_{lm}(R)=
(1/2)3^{-2/3}\pi^{-4/3}\G^{-1}(2/3)
|\sg'(w_k)|^{-1/3}
\sum_{n\in L(\nu_k)}
(-1)^{n_1}|n|^{-4/3}\sin(2\pi Y(n)R)
\eqno (5.6)
$$
is a periodic function of $R$.

Assume now that $\nu_k$ is Diophantine,
$$
|n\cdot\nu_k|\ge C|n|^{-1-\zeta},\quad 0<\zeta<1.
\eqno (5.7)
$$
In this case we put $\g=1$. Define
$$
E_{lm}^{(1)}(R;N)=\sum_{n\in \Z^2\cap L_1,\,0<|n|\le N}K(n).
$$
Let us prove that
$$
\sup_{1\le R<\infty} R^{-1/2}
|E_{lm}^{(1)}(R)-E_{lm}^{(1)}(R;N)|\le CN^{(\zeta-1)/4}.
\eqno (5.8)
$$
Indeed, $|\Ai(y)|\le C|y|^{-1/4}$ and $|a_k(\a)|\ge C_0
|\a-\a(\nu_k)|$,
hence (3.7) implies that
$$
|\tc_m(\xi)|\le C|\xi|^{-3/2}|\a(\xi)-\a(\nu_k)|^{-1/4}
$$
and
$$
|K(n)|\le CR^{1/2}|n|^{-3/2}|\a(n)-\a(\nu_k)|^{-1/4}.
$$
Therefore
$$\eqalign{
&R^{-1/2}|E_{lm}^{(1)}(R)-E_{lm}^{(1)}(R;N)|
=R^{-1/2}\left|\sum_{n\in \Z^2\cap L_1,\;|n|>N}
K(n)\right|\cr
&\le C\sum_{n\in \Z^2\cap L_1,\;|n|>N}
|n|^{-3/2}|\a(n)-\a(\nu_k)|^{-1/4}.\cr}
\eqno (5.9)
$$
Due to the Diophantine condition (5.7),
$$
|\a(n)-\a(\nu_k)|\ge C|n|^{-(2+\zeta)}.
$$
Let $J\ge N$. Order all $n\in\Z^2\cap L_1$
with $1\le|n|\le 2J$ and $n\cdot\nu_k>0$ in the increasing order
of $|\ar n\cdot\nu_k|=|n_2\nu_{k1}-n_1\nu_{k2}|$:
$$
|\ar{n^{(1)}}\cdot\nu_k|\le
|\ar{n^{(2)}}\cdot\nu_k|\le
\dots
$$
Then (5.7) implies that
$$
|\ar{n^{(1)}}\cdot\nu_k|\ge C_0J^{-(1+\zeta)}
$$
and
$$
|\ar{n^{(j+1)}}\cdot\nu_k|\ge
|\ar{n^{(j)}}\cdot\nu_k|+C_0J^{-(1+\zeta)},\quad j\ge 1,
$$
so that
$$
|\ar{n^{(j)}}\cdot\nu_k|\ge C_0jJ^{-(1+\zeta)},
$$
hence
$$
|\a(n^{(j)})-\a(\nu_k)|\ge C_1jJ^{-(2+\zeta)}
$$
and
$$
\sum_{n\in \Z^2\cap L_1,\; J\le|n|\le 2J}
|n|^{-3/2}|\a(n)-\a(\nu_k)|^{-1/4}
\le CJ^{-3/2}J^{(2+\zeta)/4}\sum_{j=1}^{C_0J}j^{-1/4}
\le C_1J^{(\zeta-1)/4}.
$$
This implies
$$
\sum_{n\in \Z^2\cap L_1,\; N>|n|}|n|^{-3/2}|\a(n)-\a(\nu_k)|^{-1/4}
\le CN^{(\zeta-1)/4},
$$
and hence (5.8) follows from (5.9). Let us evaluate now
$E_{lm}^{(2)}(R)$.

{}From (3.9)
$$\eqalign{
|E_{lm}^{(2)}(R)&-E_{lm}^+(R)-E_{lm}^-(R)|\cr
&\le CR^2\SL\pn|\tpd|\cdot|Rn|^{-5/2}
|\a(n)-\a(\nu_k)|^{-7/4}\cr}
\eqno (5.11)
$$
where
$$
E^\pm_{lm}(R)=\SL K^\pm(n)
$$
and
$$
K^\pm(n)=
R^{1/2}\PI\pn\tpd\t(n)|n|^{-3/2}|\sg_\pm(n)|^{-1/2}(-1)^{n_1}
\CSpm.\eqno (5.12)
$$

The RHS in (5.11) is estimated as follows. Let $d(\xi)
=\dist(\xi,L)$. Then
$$\eqalign{
\RHS&\le CR^{-1/2}\SL\pn|\tpd|\cdot|n|^{-5/2}
d(n)^{-7/4}|n|^{7/4}\cr
&=CR^{-1/2}\SL\pn|\tpd|\cdot|n|^{-3/4}d(n)^{-7/4}.
\cr}
$$
Since for $p=1,2,\dots,$
$$
\sum_{n\in\Z^2\setminus L_\g,\;p\le|n|\le p+1}\pn d(n)^{-7/4}
\le C,
$$
and
$$
\sup_{p\le|n|\le p+1}|\tpd|\le C(1+p\de)^{-5},
$$
we obtain
$$
\RHS\le CR^{-1/2}\sum_{p=1}^\infty (1+p\de)^{-5}p^{-3/4}
\le C_0 R^{-1/2}\de^{-1/4}=C_0 R^{-1/2}T^{1/4},
$$
so that
$$
|E_{lm}^{(2)}(R)-E_{lm}^+(R)-E_{lm}^-(R)|
\le CR^{-1/2}T^{1/4}.\eqno (5.13)
$$

Define
$$\eqalign{
F^\pm_{lm}(R)&=R^{-1/2}E_{lm}^\pm(R)=R^{-1/2}\SL K^\pm(n),\cr
F_{lm}^\pm(R;N,\de)
&=R^{-1/2}\sum_{n\in\Z^2\setminus L_\g,\;|n|\le N}
K^\pm(n)\cr}
\eqno (5.14)
$$
The central point in our proof is

{\bf Lemma 5.1.} {\it For all $N,T\ge 1$,
$$
\IO|F^\pm_{lm}(R)-F^\pm_{lm}(R;N,\de)|^2dR
\le C(N^{-1/3}+T^{-1/4}), \quad \de=T^{-1}.
$$}

We will give the proof of Lemma 5.1 below, in the end of this
section, and now let us derive Lemma 4.6 from Lemma 5.1.

Assume that $\nu_k$ is rational. Define
$$
F_{lm}(R)=R^{-1/2}(E_{lm}(R)-R^{2/3}\Phi_{lm}(R))
$$
with $\Phi_{lm}(R)$ given in (5.6). Then by (5.5)
$$
F_{lm}(R)=R^{-1/2}E_{lm}^{(2)}(R)+O(R^{-1/6}),\quad 1\le R\le T,
$$
and by (5.13)
$$
R^{-1/2}E_{lm}^{(2)}(R)=F^+_{lm}(R)+F^-_{lm}(R)+O(R^{-1}T^{1/4}),
\quad 1\le R\le T,
$$
so that
$$
F_{lm}(R)=F^+_{lm}(R)+F^-_{lm}(R)+O(R^{-1/6}+R^{-1}T^{1/4}),\quad 1\le R\le T,
$$
which implies that
$$
\liT\IO|F_{lm}(R)-F^+_{lm}(R)-F^-_{lm}(R)|^2dR=0.\eqno (5.15)
$$
Define
$$
F_{lm}^\pm(R;N)
=\sum_{n\in\Z^2\setminus L_\g,\;|n|\le N}
K_\pm(n)
$$
where
$$
K_\pm(n)=\PI
\pn\t(n)
|n|^{-3/2}|\sg_\pm(n)|^{-1/2}(-1)^{n_1}
\CSpm.
\eqno (5.16)
$$
Observe that $F^\pm_{lm}(R;N)$ is a finite trigonometric
sum, hence
$$
|F^\pm_{lm}(R;N,\de)-F^\pm_{lm}(R;N)|\le C(N)\de
=C(N)T^{-1}.
$$
Therefore Lemma 5.1 implies that
$$
\lT\IO|F^\pm_{lm}(R)-F^\pm_{lm}(R;N)|^2dR\le CN^{-1/3},
$$
hence from (5.15) we deduce that
$$
\lT\IO|F_{lm}(R)-F^+_{lm}(R;N)-F^-_{lm}(R;N)|^2dR\le CN^{-1/3}.
\eqno (5.17)
$$
This implies that $F_{lm}(R)\in B^2$ and (4.19) is the Fourier
expansion of $F_{lm}(R)$. For rational $\nu_k$ Lemma 4.6 is proved.

In the case of Diophantine $\nu_k$ we define
$F_{lm}(R)=R^{-1/2}E_{lm}(R)$. Then
$$
F_{lm}(R)=F_{lm}^{(1)}(R)+F_{lm}^{(2)}(R),\quad
F_{lm}^{(j)}(R)=R^{-1/2}E_{lm}^{(j)}(R),\quad j=1,2,
$$
and similarly to (5.17) we have that
$$
\lT\IO|F_{lm}^{(2)}(R)-F^+_{lm}(R;N)-F^-_{lm}(R;N)|^2dR\le CN^{-1/3}
\eqno (5.18)
$$
with $F^\pm_{lm}(R;N)$ defined in (5.16).

Then, (5.8) implies that
$$
\lT\IO|F^{(1)}_{lm}(R)-F_{lm}^{(1)}(R;N)|dR
\le CN^{(\zeta-1)/4},\eqno (5.19)
$$
where
$$
F^{(1)}_{lm}(R;N)=R^{3/2}\sum_{n\in L_\g,\;0<|n|\le N}
\pn\tpd\tc(2\pi nR)(-1)^{n_1}.
$$
This is a finite sum and from (3.9) we obtain that
$$
|F_{lm}^{(1)}(R;N)-\hat F^+_{lm}(R;N)-\hat F^-_{lm}(R;N)|
\le C(N)(R^{-1}+T^{-1})\eqno (5.20)
$$
with
$$
\hat F^\pm_{lm}(R)
= \sum_{n\in\Z^2\cap L_\g,\;0<|n|\le N}
K_\pm(n).
$$
It follows from (5.18)--(5.20) that
$$\eqalign{
\lT&\IO|F_{lm}(R)-
F^+_{lm}(R;N)-\hat F^+_{lm}(R;N)\cr
&-F^-_{lm}(R;N)-\hat F^-_{lm}(R;N)|^2dR\le CN^{(\zeta-1)/4},
\cr}
$$
which proves Lemma 4.6 for Diophantine $\nu_k$.

{\it Proof of Lemma 5.1.} To simplify notations we will omit
$\pm$ in sub-- and superscripts. From (5.12) and (5.14)
$$\eqalign{
\IO&|F_{lm}(R)-F_{lm}(R;N,\de)|^2dR
=\pi^{-2}\sum_{n,n'\in\Z^2\setminus L_\g;\;|n|,|n'|>N}
\pn\psi_l(n')\cr
&\tpd\tp(2\pi n'\de)
\t(n)\t(n')|n|^{-3/2}|n'|^{-3/2}|\sg(n)|^{1/2}|\sg(n')|^{1/2}
I(n,n'),\cr}
\eqno (5.21)
$$
where
$$
I(n,n')=\IO\cos(2\pi Y(n)R-\phi(n))
\cos(2\pi Y(n')R-\phi(n'))dR.
$$
Observe that
$$
|I(n,n')|\le\min\{1,T^{-1}|Y(n)-Y(n')|^{-1}\}.
\eqno (5.22)
$$
In addition,
$$
|\tpd|\le C(1+Y(n)\de)^{-5},\quad\quad C|n|\le Y(n)\le C'|n|,
$$
and the RHS of (5.21) is symmetric in $n,n'$. This implies that
$$
\IO|F_{lm}(R)-F_{lm}(R;N,\de)|^2dR
\le C_0\sum_{n,n'\in\Z^2\setminus L_\g;\; Y(n)\ge Y(n')>\b N
}H(n,n')\eqno (5.23)
$$
with
$$\eqalign{
H(n,n')&=\pn\psi_l(n')(1+Y(n)\de)^{-5}\t(n)\t(n')
Y(n)^{-3/2}Y(n')^{-3/2}\cr
&|\sg(n)|^{-1/2}|\sg(n')|^{-1/2}
\min\{1,T^{-1}(Y(n)-Y(n'))^{-1}\}\cr}
$$
and $\b>0$. First we estimate $\sum_{(n,n')\in S(N)}
H(n,n')$ with
$$
S(N)=\{ n,n'\in\Z^2\setminus L_\g;\;Y(n)\ge Y(n')>\b N;
\; Y(n)-Y(n')\ge 1\}.\eqno (5.24)
$$
Let us fix some $n$ with $Y(n)\ge\b N$ and define
the layers
$$\eqalign{
S(N,n,j)=\{n'\:n'\in\Z^2\setminus L_\g;\;
Y(n')\ge\b N;\;&j+1\ge Y(n)-Y(n')\ge j\},\cr
&1\le j\le Y(n)-\b N.\cr}
$$
Observe that $|\sg(n)|\le C|\a(n)-\a(\nu_k)|^{-1/2}$, hence
$$
\sum_{n'\in S(N,n,j)}\psi_l(n')\t(n')|\sg(n')|^{1/2}
\le C\Area S(N,n,j)\le C_0(Y(n)-j)
$$
and therefore
$$\eqalign{
&\sum_{n'\in S(N,n)}\psi_l(n')\t(n')|\sg(n')|^{1/2}
Y(n')^{-3/2}T^{-1}(Y(n)-Y(n'))^{-1}\cr
&\le C\sum_{j=1}^{Y(n)-1}(Y(n)-j)^{-1/2}T^{-1}j^{-1}
\le C_0T^{-1}Y(n)^{-1/2}\log Y(n),\cr}
$$
where $S(N,n)=\cup_j S(N,n,j)$. Making a summation in $n$
we obtain now that
$$\eqalign{
&\sum_{(n,n')\in S(N)}H(n,n')
\le CT^{-1}\sum_{n\in\Z^2\setminus L_\g;\; Y(n)>\b N}
\pn(1+Y(n)\de)^{-5}\cr
&\times Y(n)^{-2}\log Y(n)|\sg(n)|^{-1/2}.\cr}
$$

Define the layers
$$
S_j=\{ n\in\Z^2\setminus L_\g,\; j\le Y(n)\le j+1\},\quad
j>\b N.
$$
Then
$$
\sum_{n\in S_j}\pn\t(n)|\sg(n)|^{-1/2}\le C
\Area S_j\le C_0j,
$$
hence
$$\eqalign{
&\sum_{(n,n')\in S(N)}H(n,n')\le CT^{-1}
\sum_{j=\b N}^\infty (1+j\de)^{-5} j^{-1}\log j\cr
&\le C_0T^{-1}|\log\de|^2
=C_0T^{-1}\log^2T.\cr}\eqno (5.25)
$$

It remains to estimate $\sum_{(n,n')\in S_0(N)}H(n,n')$ where
$$
S_0(N)=\{n,n'\in\Z^2\setminus L_\g;\;Y(n)\ge Y(n')>\b N;\;
1\ge Y(n)-Y(n')\ge 0\}.
$$
Let us fix some $n$ with $T\ge Y(n)>\b N$. Define the layers
$$\eqalign{
S_0(N,n,j)=\{n'\in\Z^2\setminus L_\g\: Y(n')>\b N;\;
(j+1)T^{-1}&\ge Y(n)-Y(n')\ge jT^{-1}\},\cr
&j=0,1,\dots,T.\cr}
$$
To estimate the sum over $S_0(N,n,j)$ we use the following

{\bf Lemma 5.2} (2/3--estimate). {\it Let
$$
Y(\xi)=|\xi| f(|\a(\xi)-\a_0|^{1/2}),\eqno (5.26)
$$
where $f(t)\in C^\infty([0,\ep])$, $\ep>0$, $f(t)>0$,
$f'(0)=0$,
$f'''(0)\not=0$. Then if $\psi(\a)\in C^\infty([\a-0,\a_0+\ep])$
with $\supp \psi(\a)$ near $\a_0$ then
$$
\sum_{n\in \Pi(R),\;\dist(n,L)>1}
\psi(\a(n))|\a(n)-\a_0|^{-1/2}\le CR^{2/3},\eqno (5.27)
$$
where
$$
\Pi(R)=\{ n\in\Z^2\:\a_0\le\a(n)\le\a_0+\ep;\;
R\le Y(n)\le R+R^{-1/3}\}
$$
and $L=\{\xi\:\a(\xi)=\a_0\}$.}

{\it Remark.} If we step away from $\a_0$ putting
$\a_0+\ep_0\le\a(n)\le\a_0+\ep$,
$\ep_0>0$, in $\Pi(R)$ (instead of $\a_0\le\a(n)\le\a_0+\ep$), then
(5.27) reduces to a well--known 2/3--estimate (see, e.g., [CdV1]).

Proof of Lemma 5.2 is given in Appendix to the paper.

With the help of Lemma 5.2 we obtain that if $T^3\ge Y(n)$ then
$$
\sum_{n'\in S_0(N,n,j)}\psi_l(n')\t(n')|\sg(n')|^{-1/2}
\le CY(n)^{2/3}
$$
(observe that $|\sg(n')|^{-1/2}\le C|\a(n')-\a(\nu_k)|^{-1/4}
\le C|\a(n')-\a(\nu_k)|^{-1/2})$, and
$$\eqalign{
&\sum_{n'\in S_0(N,n,j)}\psi_l(n')\t(n')|\sg(n')|^{-1/2}
Y(n')^{-3/2}\min\{1,T^{-1}(Y(n)-Y(n'))^{-1}\}\cr
&\le CY(n)^{-3/2}Y(n)^{2/3}T^{-1}j^{-1}Y(n)
=CY(n)^{1/6}T^{-1}j^{-1},\quad \text{if}\quad T\ge j\ge 2,
\cr}
$$
and $\le CY(n)^{-5/6}$ if $j=0,1$. Hence
$$\eqalign{
&\sum_{n'\in S_0(N,n)}\psi_l(n')\t(n')|\sg(n')|^{-1/2}
Y(n')^{-3/2}\min\{1,T^{-1}(Y(n)-Y(n'))^{-1}\}\cr
&\le C(Y(n)^{1/6}T^{-1}\log T+Y(n)^{-5/6}),\cr}
$$
where
$$
S_0(N,n)
=\{ n'\in \Z^2\setminus L_\g\: Y(n')>\b N;\;
1\ge Y(n)-Y(n')\ge 0\}.
$$
Making a summation over $n$, we
obtain
$$\eqalign{
&\sum_{(n,n')\in S_0(N);\; Y(n)\le T^3}H(n,n')
\le C\sum_{n\in\Z^2\setminus L_\g;\; T^3\ge Y(n)\ge\b N}
\pn(1+Y(n)\de)^{-5}\cr
&\t(n)Y(n)^{-3/2}|\sg(n)|^{-1/2}
(Y(n)^{1/6}T^{-1}\log Y(n)+Y(n)^{-5/6})\equiv CI_0.\cr}
$$
Consider the layers
$$
S_j=\{ n\in\Z^2\setminus L_\g,\;j+1\ge Y(n)\ge j\},\quad T^3\ge j\ge \b N.
$$
The width of $S_j$ is of order of 1, and a simple argument shows
that
$$
\sum_{n\in S_j}\pn\t(n)|\sg(n)|^{-1/2}\le Cj.
$$
This implies
$$\eqalign{
&I_0\le C\sum_{j=\b N}^{T^3}(1+j\de)^{-5}j^{-3/2}j(j^{1/6}
T^{-1}\log j+j^{-5/6})\cr
&\le C_0(\de^{-2/3}|\log\de|T^{-1}+N^{-1/3})
=C_0(T^{-1/3}\log T+N^{-1/3}).
\cr}
$$
Thus
$$
\sum_{(n,n')\in S_0(N);\; Y(n)\le T^3}H(n,n')
\le C(T^{-1/3}\log T+N^{-1/3})\eqno (5.28)
$$
The final step is to estimate
$$
\sum_{(n,n')\in S_0(N);\; Y(n)\ge T^3}H(n,n')
$$
and this is quite simple. Since for $(n,n')\in S_0(N)$,
$$
H(n,n')\le C(1+Y(n)\de)^{-5}Y(n)^{-3}|\sg(n)|^{-1/2}
|\sg(n')|^{-1/2}
$$
and
$$
\sum_{n'\: Y(n)\ge Y(n')\ge Y(n)-1;\; n'\not\in L_\g}
|\sg(n')|^{-1/2}\le CY(n),
$$
we obtain
$$\eqalign{
&\sum_{(n,n')\in S_0(N);\; Y(n)\ge T^3}H(n,n')\cr
&\le C\sum_{(n,n')\in S_0(N);\; Y(n)\ge T^3}
(1+Y(n)\de)^{-5}Y(n)^{-2}|\sg(n)|^{-1/2}\cr
&\le C_0\sum_{j=T^3}^\infty (1+j\de)^{-5}j\le C_1T^{-2}.\cr}
\eqno (5.29)
$$
{}From (5.23), (5.25), (5.28) and (5.29) Lemma 5.1
follows.

{\it Proof of Lemma 4.7.} The proof of Lemma 4.7 is
similar in main steps to the proof of Lemma 4.6.
Assume $\nu_0=n_\G(z_0)
\in\G_j^{(l)}\subset \G_m$. By (4.15)
$$
F_{lm}(R)=\sz K(n).
\eqno (5.30)
$$
with
$$
K(n)=R^{3/2}\tpd\pn\tc_m(2\pi Rn)(-1)^{n_1}
$$
Define
$$
\eqalign{
&L=\{\la \nu_0,\quad\la\in \R\},\quad
L_1=\{n\in\Z^2\:\dist (n,L)\le 1\},\cr
&F_{lm}^{(1)}(R)=
\sum_{n\in L_1,\; n\not=0}K(n),\quad\quad
F_{lm}^{(2)}(R)=
\sum_{n\not\in L_1}K(n),\cr}
$$
so that
$F_{lm}(R)=F_{lm}^{(1)}(R)+F_{lm}^{(2)}(R)$.
By (3.15)
$$
|K(n)|\le C|n|^{-3/2},
$$
hence
$$
\sum_{n\in L_1,\;|n|\ge N}|K(n)|\le CN^{-1/2}.
\eqno (5.31)
$$
Define
$$
K_0(n)=\PI|n|^{-3/2}\tpd|\sg_0(n)|^{-1/2}\t(n)\CS.
\eqno (5.32)
$$
By (3.18)
$$\eqalign{
&\sum_{n\not\in L_1,\;|n|\ge N}
|K(n)-K_0(n)|\cr
&\le CR^{-1/2}
\sum_{n\not\in L_1,\;|n|\ge N}
|\tpd|\cdot|n|^{-2}|\a(n)-\a_0|^{-1}
\le R^{-1/2}\log^2T.
\cr}
\eqno (5.33)
$$
The following lemma holds:

{\bf Lemma 5.3}. {\it For all $N,T\ge 1$,
$$
\IT\left|
\sum_{n\not\in L_1,\;|n|\ge N}
K_0(n)\right|^2dR\le C(N^{-1/3}+T^{-1/4}).
$$}

The omit the proof of this lemma since it
basically the same as the proof of Lemma 3.3 in [Ble1]
(see also the proof of Lemma 5.1 above, where a similar statement was proved in
a more complicated situation).

Lemma 3.3 implies that for a fixed $n\in\Z^2\setminus\{0\}$,
$K(n)$ converges to $K_0(n)$ as $R\to\infty$, hence
$$
\liT\IT\left|\sum_{0<|n|<N}(K(n)-K_0(n))\right|^2dR=0,
$$
and by Lemma 5.3,
$$
\liT\IT\left|F_{lm}(R)-\sum_{0<|n|<N}K_0(n)\right|^2dR=0.
$$
This proves Lemma 4.7.

\beginsection 6. Energy Levels and Closed Geodesics \par

In this section we prove Theorem 1.2. Let
$$
N(R)=\#\{E_n\le R^2\},\quad
N_{\BS}(R)=\#\{ n\in\Z^2\: Z_2(n_1+(1/2),n_2)\le R^2,\;
|n_2|\le n_1\}\eqno (6.1)
$$
and
$$
\G=\{ p\in\R^2\: Z_2(p)=1,\;|p_2|\le p_1\}.\eqno (6.2)
$$
Recall that DH is defined by (1.12) and Hypothesis D by
(4.5), (4.6).

{\bf Lemma 6.1.} {\it DH implies Hypothesis D for $\G$.}

{\it Proof.} By (2.4) $\G$ is the graph of the function
$$
p_1=g(p_2)\equiv
|p_2|+\PI\int_a^b(1-p_2^2f^{-2}(s))^{1/2}ds,\quad
f(a)=f(b)=|p_2|,\quad |p_2|\le f_{\max}.
\eqno (6.3)
$$
$g(p_2)$ is an even $C^\infty$ function, so we will assume $p_2\ge 0$.
The function
$$
{dp_1\over dp_2}=1-\PI\int_a^b
p_2f^{-2}(s)(1-p_2^2f^{-2}(s))^{-1/2}ds
\eqno (6.4)
$$
has a nice geometric interpretation.

{\bf Proposition 6.2.}
$$
{dp_1\over dp_2}\Biggr|_{p_2=I}=-\om(I),\quad I\ge 0.\eqno (6.5)
$$

{\it Proof.} An equation of $\g(I)$ is
$$
{d\f\over dl}=If^{-2}(s),\quad {ds\over dl}
=(1-I^2f^{-2}(s))^{1/2}, \eqno (6.6)
$$
where $l$ is the normal coordinate on $\g$. Hence
$$
{d\f\over ds}=If^{-2}(s)(1-I^2f^{-2}(s))^{-1/2}
$$
and
$$
\om(I)=\PI\int_a^b{d\f\over ds}ds -1
=\PI\int_a^bIf^{-2}(s)(1-I^2f^{-2}(s))^{-1/2}ds-1.
$$
Comparing this with (6.4) we obtain (6.5). Proposition
6.2 is proved.

Observe that an inflection point on $\G$ is characterized
by $\di{d^2p_1\over dp_2^2}=0$. By (6.5) this is
equivalent to $\om'(I)=0$. Similarly, the nondegeneracy
of the inflection point is characterized by $\om''(I)\not=0$.
In addition, the Diophantine condition (4.6) is equivalent to (1.12),
hence DH implies Hypothesis D. Lemma 6.1 is proved.

Lemma 6.1 implies that Theorem 4.1 holds for $N_{\BS}(R)$.
Our goal now is to find geometric interpretation of
frequencies and amplitudes in formulas (4.11) (4.12).

Consider a finite geodesic $\g$ which starts at $x_0=(s_{\max},0)$
at some angle $-\pi/2\le\a\le\pi/2$ to the direction to the north.
$\g$ is uniquely determined by $I=\sin\a$ and $l=|\g|$,
$\g=\g(I,l)$. Let $G$ be the set of all $\g(I,l)$,
$-1\le I\le 1$, $l>0$.

Assume that $\G$ is defined as in (6.3). Define two maps:
$$
p\: G\to\G,\quad\xi\: G\to\R^2\setminus\{0\},\eqno (6.6)
$$
where
$$
p\:\g=\g(I,l)\to p(\g)=(g(I),I),
\quad
\xi\:\g=\g(I,l)\to\xi(\g)=
(l\tau^{-1}(I),\om(I)l\tau^{-1}(I)).\eqno (6.7)
$$

{\bf Proposition 6.3.} {\it $p(\g)$ and $\xi(\g)$ satisfy
$$
n_\G(p(\g))=|\xi(\g)|^{-1}\xi(\g)\eqno (6.8)
$$
and
$$
p(\g)\cdot\xi(\g)=(2\pi)^{-1}|\g|.\eqno (6.9)
$$}

{\it Proof.} (6.5) implies that $n_\G(p(\g))$ is collinear
with the vector $(1,\om(I))$ as well as $\xi(\g)$, hence
(6.8) follows. To prove (6.9) observe that the both
sides of (6.9) depend linearly on $|\g|$, so it is
sufficient to prove (6.9) in the particular case
when $|\g|=\tau(I)$. In this case (6.9) reduces to
$$
g(I)+I\tau(I)=(2\pi)^{-1}\tau(I).\eqno (6.10)
$$
Since
$$
g(I)=I+\PI\int_a^b(1-I^2f^{-2}(s))^{1/2}ds,\quad
\om(I)=\PI\int_a^b If^{-2}(s)(1-I^2f^{-2}(s))^{-1/2}ds-1,
$$
and by (6.6)
$$\tau(I)=2\int_a^b{dl\over ds}ds
=2\int_a^b(1-I^2f^{-2}(s))^{-1/2}ds,
$$
(6.10) follows. Proposition 6.3 is proved.

(6.8) implies that
$$
\xi(\g)\in L_+(n_\G(p(\g))=\{\la n_\G(p(\g)),\;\la>0\}.
$$
Hence we can define the map $\pi\: G\to N_+\G$, where
$$
N_+\G=\cup_{p\in\G} L_+(n_\G(p)),
$$
as $\pi\:\g\to (p(\g),\xi(\g))$. Observe that $\pi$ is one--to--one.

{\bf Proposition 6.4.} {\it $\g\in G$ is a closed geodesic iff
$\xi(\g)\in\Z^2$. In this case
$$\eqalignno{
Y(p(\g))|\xi(\g)|&=(2\pi)^{-1}|\g|,&(6.11)\cr
|\sg(p(\g))|^{-1/2}|\xi(\g)|^{-3/2}&=|\om'(I)|^{-1/2}
\tau(I)^{3/2}|\g|^{-3/2},&(6.12)\cr
\sgn\sg(p(\g))&=\sgn\om'(I),\quad \g=\g(I,l).&(6.13)}
$$}

{\it Proof.} Observe that $\g(I,l)$, $I\not=0$,
 is a closed geodesic with $n_1$
revolutions around the axis and $n_2$ oscillations along
the meridian iff $l=|\g|=n_2\tau(I)$ and $\om(I)=(n_1/n_2)-1$, so that
$\xi(\g)=(n_2,n_1-n_2)\in\Z^2$. Similarly, $\g=\g(0,l)$ is a
closed geodesic iff $l=|\g|=n_2\tau(I)$, so that
$\xi(\g)=(n_2,0)$. This proves the first part of Proposition 6.4.

To prove (6.11) let us notice that $Y(p(\g))|\xi(\g)|
=p(\g)\cdot\xi(\g)$, hence (6.11) follows from (6.9).
Let us prove (6.12). We have:
$$
\sg(I)=-{g''(I)\over(1+(g'(I))^2)^{3/2}},
$$
and since $g'(I)=-\om(I)$,
$$
\sg(p(\g))={\om'(I)\over(1+\om(I)^2)^{3/2}}.\eqno (6.14)
$$
On the other hand, by (6.7)
$$
|\xi(\g)|=|\g|\tau(I)^{-1}(1+\om(I)^2)^{1/2},
$$
hence
$$\eqalign{
|\sg(p(\g))|^{-1/2}|\xi(\g)|^{-3/2}
&=|\om'(I)|^{-1/2}(1+\om(I)^2)^{3/4}|\g|^{-3/2}
\tau(I)^{3/2}(1+\om(I)^2)^{-3/4}\cr
&=|\om'(I)|^{-1/2}\tau(I)^{3/2}|\g|^{-3/2}.\cr}
$$
(6.12) is proved. (6.13) follows from (6.14).
Proposition 6.4 is proved.

{\it Proof of Theorem 1.2.} From Lemma 6.1 and Theorem 4.1 we obtain
that
$$
N_{\BS}(R)=AR^2+R^{2/3}\sum_{k\: n_\G(w_k)\;\text{is rational}}
\Phi_k(R)+R^{1/2}F(R),
$$
where $\Phi_k(R)$ are periodic continuous functions and
$F(R)\in B^2$. In addition, the Fourier series of $\Phi_k(R)$
and $F(R)$ are given in formulas (4.11) and (4.12),
respectively. From Theorem 2.1 we obtain now that
$$
N(R)=AR^2+R^{2/3}\sum_{k\: n_\G(w_k)\;\text{is rational}}
\Phi_k(R)+R^{1/2}\hat F(R),
$$
where
$$
\liT\IO|F(R)-\hat F(R)|^2dR=0.
$$
This implies that $\hat F(R)\in B^2$ as well and the Fourier
series of $F(R)$ and $\hat F(R)$ coincide. If we substitute
formulas of Lemma 6.4 into (4.11), (4.12) we obtain
(1.14), (1.15). Theorem 1.2 is proved.

\beginsection  Appendix. Proof of Lemma 5.2 \par

For the sake of definiteness we will assume that $f'''(0)>0$.
Let $\G^*=\{ Y(\xi)=1\}\cap\{\a_0\le\a(\xi)\le\a_0+\ep\}$ and
$\xi^0=\G^*\cap\{\a(\xi)=\a_0\}$.
Consider a basis $e_1,e_2$ on the plane such that $e_1=\xi^0$
and $e_2$ is parallel to the tangent vector to $\G^*$
at $\xi_0$ (see Fig.6). In the basis $e_1,e_2$ the
equation of $\G^*$ has the form $\xi_1=h(\xi_2^{1/2})$,
where $h(t)\in C^\infty,\; h(0)=1,\;
h'(0)=h''(0)=0,\;
h'''(0)>0$. Let us choose the length of $e_2$ in such a way
that $h'''(0)=2$, so that $h(t)=1+t^3/3+\dots$ for small
$t$. Let $e_1^\perp,e_2^\perp$ be a dual basis to $e_1,e_2$.
Let $0\le\la(t)\le 1$ be a $C^\infty$ function on a line which
is equal to $0$ near 0, and which is equal to 1 when $t>1$.

We will show that
$$
I(R)\equiv\sum_{\Z^2\cap\Pi(R)}\la(\np)\psi(\np)(\np)^{-1/2}
\le CR^{1/6},\eqno (A.1)
$$
where $n_j^\perp=n\cdot e_j^\perp,\; j=1,2$.
Observe that $|\a(\xi)-\a_0|^{-1/2}\le C_0|\xp |^{-1/2}
R^{1/2}$ when $\xi\in\Pi(R)$, hence from (A.1) Lemma 5.2
follows.

Let $Q(\xi;R)=\la(\xp)\psi(R^{-1}\xp)|\xp|^{-1/2},$
and $Q_0(\xi;R)=Q(\xi;R)$ if $\xi\in\Pi(R)$ and
$Q_0(\xi;R)=0$ otherwise, so that
$$
I(R)=\sZ Q_0(n;R).
$$
Let
$$
Q_\de(\xi;R)=\int_{\Pi(R)} Q(\eta;R)\de^{-2}\f(\de^{-1}
(\xi-\eta)) d\eta
=Q_0(\cdot;R)*(\de^{-r2}\f(\de^{-1}\cdot))(\xi),\quad
\de>0,
$$
where $\f(\xi)\in C_0^\infty(\R^2),\;\f(\xi)\ge 0,\;
\f(\xi)=0$ if $|\xi|>1$ and $\int_{\R^2}\f(\xi) d\xi=1$, and
$$
I_\de(R)=\sZ Q_\de(n;R).
$$
Put
$$
\de=R^{-1/3}.\eqno (A.2)
$$
Observe that
$$
Q_0(\xi;R)\le CQ_\de(\xi;R),
$$
so $I(R)\le CI_\de(R)$ and to prove (A.1) it is sufficient
to show that
$$
I_\de(R)\le CR^{1/6}.\eqno (A.3)
$$

By the Poisson summation formula
$$
I_\de(R)=\sZ\tpd\tQ(2\pi n;R).\eqno (A.4).
$$
The term $n=0$ is
$$
\tQ(0;R)=\int_{\Pi(R)}Q(\xi;R) d\xi\le CR^{-1/3}
\int_0^{\ep R}|\xp|^{-1/2}d\xp\le C_0R^{1/6},
$$
hence we may consider only $n\not=0$. Let
$p_1=p\cdot e_1,\;p_2=p\cdot e_2$.
If $p_2\ge\g|p_1|,\;\g>0$, then
$$
\tQ(p;R)=\int_{\Pi(R)}Q(\xi;R)\exp(ip\xi)d\xi
\eqno (A.5)
$$
can be estimated in the following way.

First we integrate in (A.5) by the lines $|p|^{-1}p\cdot\xi=c$
and then in $c$, so that
$$
\tQ(p;R)=\II\exp(i|p|c)S(c;R)dc,
$$
where
$$
S(c;R)=\int_{\Pi(R)\cap\{|p|^{-1}p|\xi|=c\}}
\la(\xp)\psi(R^{-1}\xp)|\xp|^{-1/2}d\xi.
$$
If $\ep\ll\g$, then the lines $|p|^{-1}p\cdot\xi=c$
cross
$\Pi(R)$ transversally, which implies that
$$
S(c;R)=R^{-1/3}\la_0(c-c_0)\psi_0(R^{-1}(c-c_0);R)|c-c_0|^{-1/2},
$$
where $c_0=R|p|^{-1}p\cdot\xi^0$, $\la_0(t)\in C^\infty,\;
\la_0(t)=0$ in vicinity of 0, and $\la_0(t)=1$ in
vicinity of $\infty$ and $\psi_0(t;R)\in C_0^\infty
([0,\infty))$ has a limit in $C_0^\infty$--topology
as $R\to\infty$. Therefore
$$
|\tQ(p;R)|\le C R^{-1/3}|p|^{-5}
$$
and
$$
\left|\sum_{n\:|n\cdot e_2|\ge\g|n\cdot e_1|,\;n\not=0}
\tQ(2\pi n;R)\right|\le CR^{-1/3}.\eqno (A.6)
$$

The main difficulty is to estimate $\tQ(p;R)$ when $|p_2|
<\g p_1,\g>0$. We have:
$$\eqalign{
\tQ(p;R)
&=J\int_0^\infty dt\exp(ip_2t)\la(t)\psi(R^{-1}t)
t^{-1/2}\int_a^b ds\exp(ip_1s),\cr
&t=\xp,\;\; s=\xi_1^\perp,
\;\; a=h(t^{1/2};R),\;\; b=h(t;R+R^{-1/3}),\cr}
\eqno (A.7)
$$
where $J$ is Jacobian and $h(t;R)=Rh(R^{-1}t)$.
Let
$$
\tQ_0(p;R)
=J\int_0^\infty dt\exp(ip_2t)\psi(R^{-1}t)
t^{-1/2}\int_a^b ds\exp(ip_1s),\eqno (A.8)
$$
The difference $\tQ_1(p;R)=\tQ_0(p;R) -\tQ(p;R)$ can be estimated as
follows.

Observe that
$$
\left|\int_0^\infty dt\exp(ip_2t)(1-\la(t))\psi(R^{-1}t)t^{-1/2}
\right|\le C(1+p_2)^{-1/2},
$$
hence
$$
|\tQ_1(p:R)|\le C(1+p_2)^{-1/2}R^{-1/3}
$$
and
$$\eqalign{
&\left|\sum_{n\:|n\cdot e_2|<\g|n\cdot e_1|,\;n\not=0}
\tQ_1(2\pi n;R)\right|\le CR^{-1/3}\sZ|\tpd|
(1+|n\cdot e_2|)^{-1/2}\cr
&\le C_0R^{-1/3}\de^{-3/2}
=C_0R^{1/6}.\cr}
$$
Thus it remains to estimate a similar sum with $\tQ_0
(2\pi n;R)$.

Let us integrate in $s$ in (A.8) and make the change of
variable $u=(R^{-1}t)^{1/2}$. This gives
$$
\tQ_0(p;R)=2J(ip_1)^{-1}R^{1/2}\bigl(W(p;R)-U(p;R)\bigr),
$$
where
$$
\eqalign{
&U(p;R)=\int_0^\infty du\exp(iRp_2u^2)\psi(u^2)
\exp(iRp_1h(u)),\cr
&W(p;R)=\int_0^\infty du\exp(iRp_2u^2)\psi(u^2)
\exp[iR(1+R^{-4/3})p_1h((1+R^{-4/3})u)].\cr}
\eqno (A.9)
$$
Let us evaluate first
$$
U(p;R)=\int_0^\infty du\exp[iRp_1(yu^2+h(u))]\psi(u^2),
\quad y=p_2/p_1.\eqno (A.10)
$$
There exists a $C^\infty$ change of variable
$T=T(u,y)$ such that $T(0)=0$,
$\di{\partial T\over\partial u}(0)={\partial T\over
\partial y}(0)=1$ and
$$
yu^2+h(u)=b(y)+a(y)T+T^3/3,\eqno (A.11)
$$
where $a(y),b(y)\in C^\infty,\; a(0)=b(0)=0$
(see [H\"or2]). In addition,
$$
a(y)+T^2(0,y)=0,
$$
which follows if we differentiate both
sides of (A.11) at $u=0$. After this change of variable
$U(p;R)$ reduces to
$$
U(p;R)=\exp(iRp_1b(y))\int_{c(y)}^\infty
\exp(iRp_1(-c^2(y)T+T^3/3))\psi_0(T,y){\partial u\over
\partial T}(T,y)dT,
\eqno (A.12)
$$
with $c(y)=T(0,y)$ and $\psi_0(T,y)=\psi(t^2(T,y))$.
Following [H\"or2] let us divide $(\partial u/\partial T)(T,y)$
by $-c^2(y)+T^2$ with a remainder:
$$
{\partial u\over\partial T}(T,y)=r(T,y)(-c^2(y)+T^2)+r_0(y)
+r_1(y)T.
$$
Then the first term after substitution into (A.12) allows
integration by parts, which gives an extra $(Rp_1)^{-1}$,
and the other two terms give the main contribution to
$U(p;R)$:
$$\eqalign{
U(p;R)&=\exp(iRp_1b(y))\Bigl\{(Rp_1)^{-1/3}
V(c(y)(Rp_1)^{1/3})r_0(y)
+(Rp_1)^{-2/3}V_0(c(y)(Rp_1)^{1/3})r_1(y)\Bigr\}\cr
&+O((Rp_1)^{-1}),\cr}
\eqno (A.13)
$$
where
$$
V(x)=\int_x^\infty\exp[i(-x^2T+T^3/3)]dT,
\quad V_0(x)=\int_x^\infty T\exp[i(-x^2T+T^3/3)]dT.
$$
The method of stationary phase gives the asymptotics
of $V(x),V_0(x)$ when $|x|\to\infty$, and from this
asymptotics we obtain
$|V(x)|\le C(1+|x|)^{-1/2}$ and $|V_0(x)|\le C(1+|x|)^{1/2}$.
Therefore
$$
|U(p;R)|\le C|Rp_1|^{-1/3}\min\{1,(|p_2/p_1|\,|Rp_1|^{1/3})^{-1/2}\}
=C\min\{R^{-1/2}|p_2|^{-1/2},R^{-1/3}|p_1|^{-1/3}\}.
$$
A similar estimate holds for $W(p;R)$ and finally we obtain
$$
|\tQ_0(p;R)|\le C|p_1|^{-1}
\min\{|p_2|^{-1/2},R^{1/6}|p_1|^{-1/3}\}.
$$
Hence
$$
\left|\sum_{|n\cdot e_2|\le 1}\tpd\tQ_0(2\pi n:R)\right|
\le CR^{1/6}\sum_{|n\cdot e_2|\le 1}|n\cdot e_1|^{-4/3}\le
C_0R^{1/6},
$$
and
$$\eqalign{
&\left|\sum_{1<|n\cdot e_2|\le \g|n\cdot e_1|}\tpd\tQ_0(2\pi n:R)\right|
\le C\sum_{1<|n\cdot e_2|\le\g|n\cdot e_1|}|\tpd|
|n\cdot e_1|^{-1}|n\cdot e_2|^{-1/2}\cr
&\le C_0\de^{-1/2}=C_0R^{1/6}.\cr}
$$
Lemma 5.2 is proved.

\vskip 3mm

{\it Acknowledgements.} The author thanks Freeman Dyson, Dennis
Hejhal and Peter Sarnak for useful discussions of the paper.
This work was done at the Institute for Advanced Study,
Princeton,
and the author is grateful to IAS for financial support
during his stay at the Institute. The work was also supported
by the Ambrose Monell Foundation.

\beginsection References \par

\item {[B\'er]} P. H. B\'erard, On the wave equation on a compact Riemannian
manifold without conjugate points, {\it Math. Z.}, {\bf 155},
249--276 (1977).

\item {[Bes]} A. S. Besicovitch, {\it Almost periodic functions}, Dover
Publications, New York, 1958.

\item {[Ble1]} {P. M. Bleher, On the distribution of the number of
lattice points inside a family of convex ovals, {\it Duke Math. Journ.}
{\bf 67}, 3, 461--481 (1992).}

\item {[Ble2]} {P. M. Bleher,
Distribution of the error term in the Weyl asymptotics for the Laplace
operator on a two--dimensional torus and related lattice problems,
Preprint, Institute for Advanced Study, IASSNS-HEP-92/80, 1992
(to appear in {\it Duke Math. Journ.}).}

\item {[BCDL]} P. M. Bleher, Zh. Cheng, F. J. Dyson and J. L. Lebowitz,
Distribution of the error term for the number of lattice points
inside a shifted circle, Preprint, Inst. Adv. Study,
IASSNS-HEP-92/10, 1992 (to appear in {\it Commun. Math. Phys.}).

\item {[BL]} P. M. Bleher and J. L. Lebowitz,
Energy--level statistics of model quantum systems:
universality and scaling in a lattice--point problem,
Preprint, Inst. Adv. Study, IASSNS-HEP-93/14, 1993.

\item {[CdV1]} {Y. Colin de Verdi\`ere, Nombre de points entiers
dans une famille homoth\`etique de domaines de $\bold R^n$,
{\it Ann. Scient. \'Ec. Norm. Sup.}, $4^e$ s\'erie,
{\bf 10}, 559--576 (1977).}

\item {[CdV2]} Y. Colin de Verdi\`ere, Spectre conjoint
d'op\'erateurs pseudo--diff\'erentiels qui commutent. II.
Le cas int\'egrable,
{\it Math. Z.,} {\bf 171}, 51--73 (1980).

\item {[DG]} J. J. Duistermaat and V. W. Guillemin,
The spectrum of positive elliptic operators and
periodic bicharacteristics, {\it Inventiones math.},
{\bf 29}, 39--79 (1975).

\item {[Har1]} G. H. Hardy, The average order of the arithmetical
functions $P(x)$ and $\Delta(x)$, {\it Proc. London Math. Soc.,}
{\bf 15}, 192--213 (1916).

\item {[Har2]} G. H. Hardy, On Dirichlet's divisor problem,
{\it Proc. London Math. Soc.,}
{\bf 15}, 1--25 (1916).

\item {[HR]} D. A. Hejhal and B. Rackner,
On the topography of Maass waveforms for PSL(2,$\bold Z$):
experiments and heuristics, Preprint, University of
Minnesota,
UMSI 92/162, 1992.

\item {[Hla]} E. Hlawka, \"Uber Integrale auf konvexen K\"orpern.
I, {\it Monatsh. Math.}, {\bf 54}, 1--36 (1950).

\item {[H\"or1]} L. H\"ormander, The spectral function of an
elliptic operator, {\it Acta Math.}, {\bf 121}, 193--218 (1968).

\item {[H\"or2]} L. H\"ormander, {\it The analysis of linear partial
differential operators. I. Distribution theory and Fourier
analysis.} Springer--Verlag, Berlin e.a., 1983.

\item {[Hux]} M. N. Huxley, Exponential sums and lattice points. II.
Preprint, University of Wales College of Cardiff,
Cardiff, 1992.

\item {[Lan]} E. Landau, {\it Vorlesungen \"uber Zahlentheorie},
Chelsea, New York, 1969.

\item {[LS]} W. Luo and P. Sarnak, Number variance for arithmetic
hyperbolic surfaces, Preprint, Princeton University, Princeton, 1993.

\item {[Ran]} B. Randol, A lattice--point problem, {\it Trans. Amer. Math.
Soc.}
{\bf 121}, 257--268 (1966).

\item {[Sar]} P. Sarnak, Arithmetic quantum chaos, Schur Lectures,
Tel Aviv, 1992.

\item {[Sie]} W. Sierpinski, {\it Oeuvres Choisies}, Vol. I,
P.W.N., Warsaw, 1974, 73--108.

\vfill\eject

\beginsection Figures Captions \par

Fig.1. Geodesic on a surface of revolution.

Fig.2. The phase function $\om(I)$ for
different surfaces of revolution. Cross--sections
of the surfaces of revolution are shown
in the lower part of the figure.

Fig.3. Curve $\G$ for the surfaces of revolution
shown on Fig.2.

Fig.4. Sectorial domain with points of inflection
on the boundary.

Fig.5. Local structures of $\G$.

Fig.6. Basis $e_1,e_2$.

\bye